\documentclass[11pt]{article}

\usepackage{amssymb}

\newtheorem{theorem}{Theorem}
\newtheorem{definition}[theorem]{Definition}
\newtheorem{lemma}[theorem]{Lemma}
\newtheorem{proposition}[theorem]{Proposition}
\newtheorem{theo}{Theorem}
\newtheorem{coroa}{Corollary A}
\newtheorem{corob}{Corollary B}
\newtheorem{theow}[theorem]{Theorem}
\newtheorem{remark}{Remark}        
\newtheorem{remarks}{Remarks}
\newtheorem{example}[theorem]{Example}

\def\AA{{\cal A}}
\def\ata{\alpha^{\eta}}
\def\atar{\alpha^{\eta}_{\rho}}
\def\Atar{{\cal A}^{\eta}_{\rho}}
\def\atta{\alpha_{\eta}}
\def\BB{{\cal B}}
\def\bta{\beta^{\eta}}
\def\btar{\beta^{\eta}_{\rho}}
\def\Brta{{\cal B}^{\eta}}

\def\Btar{{\cal B}^{\eta}_{\rho}}
\def\CC{{\cal C}}

\def\Cta{{\cal C}^{\eta}}

\def\dd{{\bf d}}
\def\DD{{\cal D}}
\def\Dde{\tilde{\cal D}}
\def\Dta{{\cal D}^{\eta}}
\def\Dtad{\tilde{\cal D}^{\eta}}
\def\ee{\varepsilon}
\def\esp{{\mathbb{E}}}
\def\Eta{{\cal E}_{\eta}}
\def\ftar{\phi^{\eta}_{\rho}}
\def\gta{\gamma^{\eta}}
\def\gtar{\gamma^{\eta}_{\rho}}
\def\lacc{\left\{}
\def\lcr{\left[}
\def\lpa{\left(}
\def\lro{\lambda_{\rho}}
\def\lta{\lambda_{\eta}}
\def\lva{\left|}
\def\mta{\mu_{\eta}}
\def\mub{{\bar \mu}}
\def\mubde{{\widetilde{\bar \mu}}}
\def\mude{{\widetilde \mu}}
\def\NN{{\mathbb{N}}} 
\def\nub{{\bar \nu}}
\def\pb{{\mathbb{P}}}
\def\pstar{\psi^{\eta}_{\rho}}
\def\rl{{\mathbb{R}}}
\def\racc{\right\}}
\def\rcr{\right]}
\def\rpa{\right)}
\def\rta{\rho_{\eta}}
\def\rva{\right|}
\def\Saw{{\rm Saw}}
\def\sk{{\cal D}}
\def\ssc{{\cal S}}
\def\sso{{\overline \ssc}}
\def\sstar{s^{\eta}_{\rho}}
\def\Sta{{\cal S}^{\eta}}
\def\tphi{{\widetilde \phi}}
\def\tpsi{{\widetilde \psi}}
\def\Ttar{T^{\eta}_{\rho}}
\def\Un{{\bf 1}}
\def\uta{u_{\eta}}
\def\Utar{U^{\eta}_{\rho}}
\def\UU{{\bf U}}
\def\vta{v^{\eta}}
\def\vtaL{v^{\eta}_L}
\def\Vta{{\cal V}^{\eta}}

\def\Vtar{{\cal V}^{\eta}_{\rho}}
\def\vtar{v^{\eta}_{\rho}}
\def\Wp{{\cal W}_p}
\def\wta{w^{\eta}}
\def\xita{\Xi_{\eta}}
\def\xta{x^{\eta}}
\def\Xta{X^{\eta}}
\def\xtar{x^{\eta}_{\rho}}
\def\Yjd{Y^j}
\def\Zba{{\bar Z}}
\def\Zbde{{\tilde \Zba}}
\def\Zbta{\Zba^{\eta}}
\def\Zbtad{\Zbde^{\eta}}
\def\Zde{{\tilde Z}}
\def\Zjd{\Zde^j}
\def\zta{z^{\eta}}
\def\Zta{Z^{\eta}}
\def\Ztad{\Zde^{\eta}}
\def\Ztar{Z^{\eta}_{\rho}}
\def\Ztard{\Zde^{\eta}_{\rho}}
        
\newcommand{\dsty}{\displaystyle}
\newcommand{\fin}{\vspace{-0.3cm}
                  \begin{flushright}
                  \mbox{$\Box$}
                  \end{flushright}
                  \noindent}

\begin{document}

\begin{center}

\huge
{\bf Small deviations in $p$-variation for multidimensional L\'evy processes} 

\vspace{8mm}

\Large
{\bf Thomas Simon}

\end{center}

\vspace{2mm}

\begin{abstract}

\vspace{2mm}

\noindent
{Let $Z$ be an $\rl^d$-valued L\'evy process with strong finite $p$-variation for some 
$p<2$. We prove that the ''decompensated'' process $\Zde$ obtained from $Z$ by 
annihilating its generalized drift has a small deviations property in $p$-variation. This property
means that the null function belongs to the support of the law of $\Zde$ with respect to the $p$-variation distance. 
Thanks to the continuity results of T.~J.~Lyons/D.~R.~E.~Williams \cite{lyons1} \cite{will1}, this allows us to 
prove a support theorem with respect to the $p$-Skorohod distance for canonical SDE's driven by $Z$ without any assumption 
on $Z$, improving the results of H.~Kunita \cite{kuni2}. We also give a criterion ensuring the small deviation property for $Z$
itself, noticing that the characterization under the uniform distance, which we had obtained in \cite{dev}, no more holds under
the $p$-variation distance.}
\end{abstract}

\noindent 
{\bf Keywords:} L\'evy measure - Marcus equation - $p$-variation - Support.

\vspace{2mm} 

\noindent
{\bf MSC 2000:} 60G51, 60H10  

\section{Introduction}

In a series of celebrated papers \cite{lyons1} \cite{lyons2}, T.~J.~Lyons has built a general theory of rough differential equations. One of the main interests of this theory is the possibility to solve path-wise multidimensional stochastic equations whose driving paths have finite $p$-variation only for some $p > 1$. In contrast to It\^o's theory, the equations are solved through convergence of the Picard iteration scheme, so that their solutions can be viewed as continuous functionals of the driving signal, with respect to some $p$-variation distance. Lyons' papers dealt only with the continuous case, and recently D.~R.~E.~Williams \cite{will1} \cite{will2} has extended this theory to discontinuous stochastic equations, in particular when the driving noise is a L\'evy process. Williams considered equations with jumps of It\^o and Marcus type, but got continuity results only in the latter case. The continuity problem for equations of It\^o type seems namely quite difficult, even in dimension 1.

The purpose of this paper is to apply Lyons/Williams' results to the proof of a support theorem for S.D.E.'s of Marcus type driven by a multidimensional L\'evy process. Namely, this continuity property allows us to reduce the difficult part of this kind of theorem (the approximation of an element of the support by the S.D.E. with positive probability) to a control on the driving path itself. In \cite{marcus} we had already used this idea when the underlying L\'evy process is one-dimensional. In this situation the continuity property holds with respect to the local uniform norm, so that we could appeal to the small deviation results for L\'evy processes in uniform topology which we had obtained in \cite{dev}. Here our driving L\'evy process $Z$ is multidimensional, and we restrict ourselves to the case where it has finite $p$-variation for some $p<2$. In particular we only make use of the "little theorem" of Lyons/Williams, and no control on the area process is involved.

However we make no other assumption on $Z$, so that our support theorem covers a wide class of driving L\'evy processes without Gaussian part, and improves significantly the results of Kunita \cite{kuni2} on the subject. Our description of the support is also simpler, and more naturally related to the geometry of the underlying L\'evy measure $\nu$, as in Tortrat's famous article \cite{tortr1}: up to finitely many jumps and some fixed drift, this support is made out of a family of O.D.E.'s driven by functions with regular $p$-variation which are valued in the subspace of $\rl^d$ consisting in the completely asymptotic directions of $\nu$ - see (2) below for details. Notice finally that thanks to the continuity property, we obtain the support theorem in a stronger topology than the local Skorohod one, taking into account the $p$-variation, and which we call the {\em $p$-Skorohod topology.}

The core of this article consists in the proof of the small deviations property in $p$-variation norm for the ``decompensated process'' $\Zde$ obtained from $Z$ after annihilating its generalized drift. This generalized drift is just the sum of the usual drift and of the projection of finite 1-variations of the compensator in the L\'evy-Khintchine formula - in particular the projection of $\Zde$ with finite 1-variations is the sum of its jumps. The small deviation property simply means that the null function belongs to the support of $\Zde$ with respect to the $p$-variation norm. In \cite{dev} we had already obtained this property for $\Zde$ under the uniform norm. Here the proof is somewhat analogous, but much more delicate. Roughly, denoting by $L$ the above completely asymptotic subspace and by $\Pi^{}_L$ the orthogonal projection operator onto $L$, we need to approximate $\Zde$ by some ''saw-function'' whose slope is 
$$\vtaL\; =\; \Pi_L^{}\lpa\int_{\eta\leq |z|\leq 1}\!\! z\,\nu(dz)\rpa$$
for every small $\eta$ along some subsequence. In other words, we must approximate $\vtaL$ by a sum of the type
$$\sum_{i=1}^r\ata_i\xta_i$$
where $r$ is a fixed integer, $\xta_1,\ldots,\xta_r\in\mbox{Supp}\,\nu\cap\{|z|\leq\eta\}$, and $\ata_1,\ldots,\ata_r$ are minimizing integers verifying
\begin{eqnarray*}
\ata_i{\lva\xta_i\rva}^p & \rightarrow & 0
\end{eqnarray*}
for every $i = 1,\ldots, r$, as $\eta$ tends to 0 along the subsequence. The latter convergence is crucial because of the $p$-variation norm, but it is quite hard to obtain in full generality on the L\'evy measure. We overcome the difficulties with the help of elementary analysis and geometry which require a lot of care, and where strict positivity (or strict convexity) plays a central r\^ole. A useful tool is also Skorohod's absolute continuity theorem for L\'evy processes, which comes rather unexpectedly since the involved transformations are far less tractable than in the Cameron-Martin theorem.

In \cite{dev} a {\em characterization} of the small deviation property under the uniform norm for general multidimensional L\'evy processes was obtained, in terms of interactions between the drift and the projection of finite 1-variations of the L\'evy measure. Simple examples show that this characterization no more holds under the $p$-variation norm: there are L\'evy processes with finite $p$-variation which have small deviations under the uniform norm but not under the $p$-variation norm. The characterization in the case $p=1$ is easily proved to be $Z=\Zde$, that is $Z$ itself is the sum of its jumps. In the case $p>1$ and when the L\'evy measure has infinite variations in every direction, we also prove that the small deviation property in $p$-variation always holds for $Z$. In the general case when $p>1$ and the projection of $Z$ with finite $1$-variation is non trivial, we give a criterion involving the strict convexity of the asymptotic c\^one generated by Supp $\nu$, a criterion which in some sense is optimal.

The organization of this article is as follows: in Section 2 we present the framework and state the main result of this paper - the small deviation property for $\Zde$, as well as its two corollaries - the criterion mentioned above and the support theorem for Marcus equations driven by $Z$. In this section we also give some examples, which might be helpful for the understanding of the proof of the main theorem. The latter, which is unfortunately quite technical, is given in Section 4. Before that, we give in Section 3 a few lemmas concerning some deterministic functions - in particular the saw-functions which are of central use in the proof of the main result - and their $p$-variation. In Section 5 we prove the corollaries.

\section{Notations and results}

\subsection{L\'evy processes and their $p$-variation}

We work on $\rl^d$ endowed with $\lva .\rva$ any Euclidean norm. Let $p\geq 1$ and $I$ be an interval of $\rl^+$. A function $f : \rl^+ \rightarrow \rl^d$ is said to have finite (strong) $p$-variation over $I$ if 
$${\lva\lva f\rva\rva}_{I, p}\; =\; {\lpa\sup_{\tiny{t_0 < \ldots < t_k\in 
I}}\sum_{j=1}^k{\lva f(t_j) - f(t_{j-1})\rva}^p\rpa}^{1/p}\;<\; \infty.$$
If $I =[0,T]$ for some $T\geq0$ we will use the simpler notation
${\lva\lva f\rva\rva}_{T, p}$ for ${\lva\lva f\rva\rva}_{[0,T], p}$.
Notice that for every $q\geq p$,
$${\lva\lva f_0\rva\rva}_{I, \infty}\; \leq\;{\lva\lva f\rva\rva}_{I, q}\; 
\leq\;{\lva\lva f\rva\rva}_{I, p}$$
where we wrote $f_0(t) = f(t) - f(0)$ for every $t\geq 0$ and ${\lva\lva .\rva\rva}_{I, \infty}$ stands for the uniform norm over $I$. Besides ${\lva\lva .\rva\rva}_{I, p}$ is a Banach semi-norm which satisfies in particular the triangle inequality. We denote by $\Wp$ the space of functions having finite $p$-variation over every compact interval, factored by the set of constant functions. The family of semi-norms
$$\lacc {\lva\lva .\rva\rva}_{n, p},\; n\geq 1\racc$$
makes $\Wp$ into a Banach space with a norm ${\lva\lva .\rva\rva}_p$ defined in the usual way:
$${\lva\lva f\rva\rva}_p\; = \;\sum_{n \geq 1} 2^{-n}\lpa 1 \wedge {\lva\lva f\rva\rva}_{n, p}\rpa$$
for every $f\in\Wp$. A function $f : \rl^+ \rightarrow \rl^d$ is said to have {\em regular} finite $p$-variation over an interval $I$ if
$$\lim_{\ee\rightarrow 0}\lpa\sup_{\tiny{\begin{array}{c}
        t_0 < \ldots < t_k\in I\\
        |t_j - t_{j-1}|\leq \ee
        \end{array}}}\sum_{j=1}^k{\lva f(t_j) - f(t_{j-1})\rva}^p\rpa\;=\; 0.$$
Notice that this notion is only of interest for $p>1$, and that if $f$ is continuous with finite $p$-variation over $I$, then it has regular finite $q$-variation over $I$ for every $q>p$. Notice also that every function with regular finite $p$-variation is necessarily continuous.

\vspace{2mm}

We now fix $x\in\rl^d$ and denote by $\Wp(x)$ the space of functions having finite $p$-variation over every compact interval and starting from $x$. In the following we shall also work on $\rl^m$ for some $m\neq d$, and we will still denote by $\Wp(x)$ the space of functions having finite $p$-variation over every compact interval and starting from $x\in\rl^m$. 

It is well-known and easy to see that every member of $\Wp(x)$ has left and 
right limits at every point of $\rl^+$. We denote by $\sk_p(x)$ the subspace of 
$\Wp(x)$ made out of {\em c\`ad-l\`ag} functions. We endow it with the 
following distance: if $f, g\in\sk_p(x)$
$$\dd_p(f,g) \; = \;\sum_{n \geq 1} 2^{-n}\lpa 1 \wedge \dd_p^n(f,g)\rpa,$$
where for every $n\in\NN^*$ $\dd_p^n$ is defined by
$$\dd_p^n(f,g)\; =\;\inf_{\lambda\in\Lambda}\lacc\sup_{s\leq 
t}\lva\log\frac{\lambda_t -\lambda_s}{t-s}\rva + {\lva\lva k_n 
f(\lambda_.) - k_n g(.)\rva\rva}_{n+1, p}\racc,$$ 
$\Lambda$ designing the set of all continuous strictly increasing functions 
$\lambda : \rl^+ \rightarrow \rl^+$ with $\lambda_0 = 0$ and $\lambda_t\uparrow 
+\infty$ as $t\uparrow +\infty$, and $k_n$ being given by
$$k_n(t)\; =\;\lacc\begin{array}{ll}
                   1 & \mbox{if $t\leq n$}\\
                   n+1 - t & \mbox{if $n < t \leq n+1$}\\
                   0 & \mbox{if $t\geq n+1.$}
                   \end{array}
              \right.
$$
Such a $\lambda$ will be called a {\em change of time} in the sequel. Making the same considerations as in \cite{jacshy} pp. 293-294 
and using the fact that ${\lva\lva .\rva\rva}_{n+1, p}$ is a semi-norm entails that $\dd_p$ 
is actually a distance on $\Wp(x)$, which dominates the usual local Skorohod distance 
$\dd$: for every $f, g\in\sk_p(x)$,
$$\dd(f,g)\;\leq\;\dd_p(f,g).$$
In the sequel $\dd_p$ will be called the {\em $p$-Skorohod distance} and the topology induced by $\dd_p$ on $\sk_p(x)$ the {\em $p$-Skorohod topology}.

\vspace{2mm}

Let $\lacc Z_t, t \geq 0\racc$ be an $\rl^d$-valued L\'evy process starting from 
0, without Gaussian part. Its L\'evy-It\^o decomposition writes
$$Z_t\; = \; \alpha t \; + \; \int_0^t\int_{|z| \leq 1}\!\!
z\,\mude(ds,dz) \; + \; \int_0^t\int_{|z| > 1}\!\! z\,\mu(ds,dz),$$
where $\alpha\in\rl^d$, $\nu$ is a positive Borel measure on $\rl^d - \{0\}$ 
satisfying
$${\dsty \int_{\rl^d}\frac{|z|^2}{|z|^2 + 1}}\, \nu(dz) \; < \;\infty,$$ $\mu$ 
is the Poisson measure over $\rl^+\!\times\rl^d$ with intensity $ds 
\otimes\nu(dz)$, and $\mude = \mu - ds \otimes\nu$ is the compensated measure. 
Bretagnolle \cite{breta} obtained the following characterization: for every $1\leq p <2$
$$Z\in \Wp(0)\;\;\mbox{a.s.}\;\Longleftrightarrow\;\int_{|z| \leq 1}\!\! 
{|z|}^p\nu(dz)\; <\; \infty$$
(notice that the equivalence is trivial for $p= 1$). In particular every stable process has finite $p$-variation for some $p <2$. Recall that on the contrary Brownian Motion has infinite 2-variation, so that the above characterization only makes sense for L\'evy processes without Gaussian part. 

\vspace{2mm}

In the case $p > 1$, Bretagnolle got actually a sharper result: the existence of two universal constants $c_p$ and $C_p$ depending only on $p$ such that
$$c_p \int_{|z| \leq 1}\!\! {|z|}^p\nu(dz)\;\leq\;\esp\lcr {\lva\lva Z \rva\rva}^p_{1,p}\rcr\;\leq\; C_p\int_{|z| \leq 1}\!\! {|z|}^p\nu(dz)$$
if $\nu$ is concentrated on $\{|z|\leq 1\}$ and $\alpha = 0$. Using the inequality ${\lva a+b\rva}^p\;\leq\; 2^{p-1}\lpa {\lva a\rva}^p +{\lva b\rva}^p\rpa$ for any $a, b\in\rl^d$ entails easily that for every $T > 0$
$$2^{1-p}c_p T \int_{|z| \leq 1}\!\! {|z|}^p\nu(dz)\;\leq\;\esp\lcr {\lva\lva Z \rva\rva}^p_{T,p}\rcr\;\leq\; 2^{p-1}C_p T \int_{|z| \leq 1}\!\! {|z|}^p\nu(dz)$$
if $\nu$ is concentrated on $\{|z|\leq 1\}$ and $\alpha = 0$. This latter estimate yields readily the following approximation lemma - notice that the case $p=1$ is trivial:

\begin{lemma} Let $\lacc \Eta,\;\eta >0\racc$ be a family of subsets of $\rl^d$ included, for every $\eta >0$, in the ball of radius $\eta$ centered at the origin. Let $Z\in\Wp(0)$ be a L\'evy process with drift $\alpha$ and L\'evy measure $\nu$. Let $\Zta$ be the L\'evy process with same drift and L\'evy measure $\Un_{{}_{\Eta^c}}\nu$ and set $\Ztad = Z - \Zta$. Then for every $T>0$,
$$\esp\lcr{\lva\lva\Ztad\rva\rva}_{T, p}\rcr\;\longrightarrow\; 0$$
as $\eta\downarrow 0$.
\end{lemma}

The following u.c.p. lemma (also trivial in the case $p=1$) is another immediate application of Bretagnolle's estimate, thanks to the Markov inequality:

\begin{lemma} Let $\nu_0$ be a L\'evy measure on $\rl^d$ concentrated on $\{|z|\leq 1\}$, and integrating ${\lva z\rva}^p$. For every $\ee > 0$, $0 < c <1$ and $T> 0$, there exists $\eta_0 > 0$ such that for every $\eta < \eta_0$ and every L\'evy process $Z\in\Wp(0)$ with drift $\alpha = 0$ and L\'evy measure $\nu \leq \nu_0$,
$$\pb\lcr{\lva\lva\Ztad\rva\rva}_{T, p} < \ee \rcr\; >\; c$$
with the notations of Lemma 1, and where the notation $\nu \leq \nu_0$ means that $\nu (A) \leq \nu_0 (A)$ for every measurable $A\subset\rl^d$. 
\end{lemma}
  
\subsection{Small deviations in $p$-variation norm for L\'evy processes}

As in \cite{tortr1} and \cite{sharpe} we introduce the following vector 
space
$$K\; = \;\lacc x\in\rl^d / \int_{|z| \leq 1}\!\!|x*z|\,\nu(dz)\; <\;
\infty\racc,$$
where $*$ is the scalar product defining the chosen Euclidean norm on 
$\rl^d$. Notice that the vector space $L = K^{\perp}$, which can be viewed as 
the completely asymptotic direction of $\Un_{|z|\geq\eta}\nu(dz)$ as 
$\eta\downarrow 0$, depends {\em only} on $\nu$ and not on the choice of this Euclidean 
structure. In the following, every L\'evy process with the above L\'evy-It\^o 
decomposition will be said to have characteristics $(\alpha,\nu)$, and $K$ 
will always implicitly stand for the orthogonal space of $L$ with respect to 
the chosen scalar product. We define the {\em generalized drift} of a L\'evy process $Z$ with 
characteristics $(\alpha,\nu)$ by
$$\alpha^{}_{\nu}\; =\;\alpha\; -\;\int_{|z|\leq 1}\!\!\!\! z^{}_K\,\nu(dz)$$
where $z^{}_K$ is the orthogonal projection of $z$ onto $K$, so that the integral 
makes sense. We finally introduce the {\em decompensated process} $\Zde$ 
associated with $Z$, which is the L\'evy process with characteristics 
$(\alpha-\alpha^{}_{\nu},\nu)$. Equivalently
$$\Zde_t \; =\; \int_0^t\int_{|z|\leq 1}\!\!\!\! z^{}_K\,\mu(ds,dz)\; 
+\;\int_0^t\int_{|z|\leq 1}\!\!\!\! z^{}_L\,\mude(ds, dz)\; +\; \int_0^t\int_{|z| > 
1}\!\! z\,\mu(ds,dz)$$
for every $t >0$, where $z^{}_L$ denotes the orthogonal projection of $z$ onto $L$.
The main result of this paper is the following

\begin{theo} Let $Z\in\Wp(0)$ be a L\'evy process with characteristics $(\alpha,\nu)$. 
Its decompensated process $\Zde$ has the following small deviation property:
$$\pb\lcr{\lva\lva\Zde\rva\rva}_{T,p}\; < \; \ee\rcr\; > \; 0$$
for every $\ee > 0$ and $T > 0$.
\end{theo}

We now recall a few notations from \cite{tortr1}, \cite{sharpe} and \cite{dev}: for every $\eta >0$, set $\Cta$ for the closed convex c\^one with vertex $0$ generated by $\Sta = \mbox{Supp}\;\nu\;\cap\;\lacc |z| \leq \eta\racc$, and
$$\CC\; =\; \bigcap_{\eta > 0} \Cta.$$
Let $\Pi^{}_K$ be the operator of orthogonal projection onto $K$,
$$\AA^{}_K\; =\; \lpa \int_{|z|\leq 1}\!\!\!\! z^{}_K\,\nu(dz)\rpa\; - \;\overline{\Pi^{}_K \lpa\CC\rpa}\;\;\;\mbox{and}\;\;\;\BB^{}_K\; =\; \lpa \int_{|z|\leq 1}\!\!\!\! z^{}_K\,\nu(dz)\rpa\; - \;\bigcap_{\eta >0}\overline{\Pi^{}_K \lpa\Cta\rpa}.$$
It follows from the main result of \cite{dev} that if $Z\in\Wp(0)$ is a L\'evy process with characteristics $(\alpha,\nu)$, then the following equivalence holds:
$$\alpha\in\Pi^{-1}_K\lpa\BB^{}_K\rpa\;\Longleftrightarrow\;\pb\lcr{\lva\lva Z\rva\rva}_{T,\infty}\; < \; \ee\rcr\; > \; 0\;\;\;\mbox{for every $\ee > 0$ and $T > 0$.}$$
One could wonder if the same characterization holds under the $p$-variation norm, i.e. if one could replace $\infty$ by $p$ in the above right-hand side. However this is not true, as shows the following example.

\begin{example}{\em Consider on $\rl^2 = \{(z_1,z_2)\}$ endowed with the canonical basis $\lpa e_1, e_2\rpa$, the following measure
$$\nu(dz) \; =\; \Un_{\{0<z_1<{\lva z_2\rva}^{r}<c_r\}}{\lva z_2\rva}^{-2-q}\,dz$$
where $q$ and $r$ are such that $1<(1+q)/2<r<q<r+1$ (notice that $q$ can take any value strictly greater than $1$), and $c_r$ is the unique positive solution to $x^2 + x^{2/r} = 1$. Then $\nu$ is a jumping measure concentrated
on $\lacc |z| \leq 1\racc$, whose asymptotic subspaces are $K\; =\; \mbox{Vect}\{e_1\}$ and $L\; =\; \mbox{Vect}\{e_2\}.$
Besides, with the above notations, $\Pi^{-1}_K\lpa\AA^{}_K\rpa\; =\; \Pi^{-1}_K\lpa\BB^{}_K\rpa\; =\;\lacc z_1\leq c\racc$ where we set
$$c = \frac{1}{2r - q - 1}.$$
 Let $Z$ be the L\'evy process given by
$$Z_t\; =\; \int_0^t\int_{|z| \leq 1}\!\! z\,\mude(ds,dz)$$
for every $t>0$. With our notations, $\alpha\; = \;0\;\in\;\Pi^{-1}_K\lpa\BB^{}_K\rpa$. On the other hand, $Z\in\Wp(0)$ a.s. for every $p > 1 + q -r$ (notice that $0<q-r < 1)$. Set now $p\in ]1+q-r,r[$ (notice that $r > 1+q -r$) and let $Z^1$ (resp. $Z^2$) be the projection of $Z$ onto $K$ (resp. onto $L$). For every $0 < \ee < c/2$,
\begin{eqnarray*}
\lacc{\lva\lva Z\rva\rva}_{1,p}\; < \; \ee\racc & \subset & \lacc\sup_{t\leq 1}\lva\Delta Z_t\rva\; < \; \ee,\;\sum_{t\leq 1}\Delta Z^1_t\; > \; c/2\racc\\
& \subset & \lacc\sup_{t\leq 1}\lva\Delta Z_t\rva\; < \; \ee,\;\sum_{t\leq 1}{\lva\Delta Z^2_t\rva}^r\; > \; c/2\racc\\
& \subset & \lacc\sum_{t\leq 1}{\lva\Delta Z_t\rva}^p\; > \; \ee^{p-r}c/2\racc,
\end{eqnarray*}
which proves that
$$\pb\lcr{\lva\lva Z\rva\rva}_{1,p}\; < \; \ee\rcr\; = \; 0$$
as soon as $\ee < {\lpa c/2\rpa}^{1/r}$, since obviously
$${\lva\lva Z\rva\rva}_{1,p}\; \geq \; {\lpa\sum_{t\leq 1}{\lva\Delta Z_t\rva}^p\rpa}^{1/p}\;\;\mbox{a.s.}$$
}
\end{example}

\noindent
Nevertheless, our main result makes it possible to obtain the following criterion. We say that a closed convex c\^one with vertex 0 in $\rl^d$ is {\em strictly convex} if it contains no half-plane. 

\begin{coroa} Let $Z\in\Wp(0)$ be a L\'evy process on $\rl^d$ with characteristics $(\alpha,\nu)$. With the above notations,

\vspace{2mm}

\noindent
{\em (a)} If $K =\rl^d$, then $Z$ has small deviations in $1$-variation if and only if $\alpha_{\nu} =0$. 

\vspace{2mm}

\noindent
{\em (b)} If $L = \rl^d$, then $Z$ has small deviations in $p$-variation ($p > 1$).

\vspace{2mm}

\noindent
{\em (c)} If $L\neq\rl^d, \alpha\in \Pi^{-1}_K\lpa\AA^{}_K\rpa$ and $\CC$ is strictly convex, then $Z$ has small deviations in $p$-variation ($p > 1$) . 

\vspace{2mm}

\noindent
{\em (d)} If $L\neq\rl^d$ and $\alpha\notin \Pi^{-1}_K\lpa\BB^{}_K\rpa$, then $Z$ does not have small deviations in $p$-variation ($p > 1$) .

\end{coroa}

\begin{remarks}{\em (a) Even if $\CC$ is strictly convex, the condition $\alpha\in \Pi^{-1}_K\lpa\BB^{}_K\rpa$ is not sufficient, as the following example easily shows. Consider
$$\nu(dz) \; =\; \Un_{\{0<z_1<z_2^r<c_r\}}z_2^{-(2+q)}\,dz$$
on $\rl^+\times\rl^+$, with the notations of Example 3. Here,
$$\Pi^{-1}_K\lpa\AA^{}_K\rpa\; =\;\lacc z_1 = \frac{1}{2(2r-q-1)}\racc\;\;\;\mbox{and}\;\;\;\Pi^{-1}_K\lpa\BB^{}_K\rpa\; =\;\lacc z_1\leq \frac{1}{2(2r-q-1)}\racc.$$
The process $Z$ defined as in Example 3 verifies of course $\alpha\in \Pi^{-1}_K\lpa\BB^{}_K\rpa$, but does not have small deviations in $p$-variation either.

\vspace{2mm}

\noindent
(b) If $\CC$ is strictly convex, the condition $\alpha\in \Pi^{-1}_K\lpa\AA^{}_K\rpa$ is sufficient but {\em not} necessary, as shows the following example. 
Consider on $\rl^+\times\rl^+$ 
$$\nu(dz) = \nu_1(dz) + \nu_2(dz),$$
where
$$\nu_1(dz) \; =\; \Un_{\{0<z_1<z_2^r<c_r\}}\,dz\;\;\;\mbox{and}\;\;\;\nu_2(dz) \; =\;\sum_{n\geq 1}n^q\delta_{(0,n^{-1})}(dz).$$
Here $0<q<1<r$ and $c_r$ is the unique positive solution to $x^2 + x^{2/r} = 1$. Then $\nu$ is a jumping measure concentrated
on $\lacc |z| \leq 1\racc$, whose asymptotic subspaces are $K\; =\; \mbox{Vect}\{e_1\}$ and $L\; =\; \mbox{Vect}\{e_2\}.$
Besides
$$\Pi^{-1}_K\lpa\AA^{}_K\rpa\; =\;\lacc z_1 = \frac{1}{2(2r+1)}\racc\;\;\;\mbox{and}\;\;\; \Pi^{-1}_K\lpa\BB^{}_K\rpa\; =\;\lacc z_1\leq \frac{1}{2(2r+1)}\racc.$$
Let $\alpha\in\rl^2$ with $\alpha_1 < \frac{1}{2(2r+1)}$ and consider $Z$ the L\'evy process with characteristics $(\alpha,\nu)$. Clearly $Z\in\Wp(0)$ a.s. if $p > q+1$. We briefly show that $Z$ has small deviations in $p$-variation norm, even though $\alpha\notin\Pi^{-1}_K\lpa\AA^{}_K\rpa$. Set $\mu_i$ for the Poisson measure on $\rl^+\times \rl^d$ with compensator $ds\otimes\nu_i(dz)$, $i=1,2$. Introduce the compound Poisson process
$$Z^1_t\; = \;\int_0^t\int_{|z| \leq 1}\!\! z\,\mu_1(ds,dz)\; =\; \sum_{s\leq t}\Delta Z^1_s$$
for every $t>0$, and let $\{T_n,Z_n\}_{n\geq 1}$ be the sequence of its successive jumping times and sizes. Let $T_0 =0$ and $S_n = T_n -T_{n-1}$ for every $n\geq 1$. Set finally 
$$\beta\; =\; \frac{1}{2(2r+1)}\, -\, \alpha_1.$$
Fix $p> q+1$, $\ee, T > 0$. We will show that
$$\pb\lcr{\lva\lva Z\rva\rva}_{T,p}\; < \; \ee\rcr\; > \; 0$$
in using the independence of $Z^1$ and $\mu_2$. Take $\eta > 0$ such that $6cT\eta^{p-r} <\ee$. Consider $P_{\eta}=(\eta^r, \eta)\,\in\,\mbox{Supp}\,\nu_1$ and $t_{\eta} = \eta^r/cT$. For every $\lambda > 0$ the event
$$\Omega_{\lambda}\; =\; \lacc \lva S_n - t_{\eta}\rva < \lambda,\; \lva Z_n - P_{\eta}\rva < \lambda,\;\;\forall\, n\, = 1,\ldots,\lpa \mbox{Ent}\lcr cT/\eta^r\rcr +1\rpa\racc$$
has positive probability. On the other hand, it follows from Lemma 11 below that 
$$\Omega_{\lambda}\;\subset\;\lacc {\lva\lva Z^1_{\eta}\rva\rva}_{T,p}\; < \; \ee/2\racc$$ 
if $\lambda$ is small enough, where we set
$$Z^1_{\eta}(t)\; = \;\sum_{s\leq t}\Delta Z^1_s\; -\; tc\lpa e_1 + \eta^{1-r}  e_2\rpa$$
for every $t >0$. But by Corollary A,
$$\pb\lcr{\lva\lva Z^2_{\eta}\rva\rva}_{T,p}\; < \; \ee/2\rcr\; > \; 0$$
where we introduced the process
$$Z^2_{\eta}(t)\; = \;\int_0^t\int_{|z| \leq 1}\!\! z\,\mude_2(ds,dz)\; +\; t\lpa \eta^{1-r}c +\alpha_2\rpa e_2$$
for every $t > 0$. Since $Z = Z^1_{\eta} + Z^2_{\eta}$ for every $\eta > 0$, with $Z^1_{\eta}$ and $Z^2_{\eta}$ independent, we finally get
$$\pb\lcr{\lva\lva Z\rva\rva}_{T,p}\; < \; \ee\rcr\; > \; 0$$
by the triangle inequality. We stress finally that the above ''approximation event'' $\Omega_{\lambda}$ will be introduced repeatedly during the proof of our main result, in various forms.}

\end{remarks}

We next give two more classical examples which fall into the scope of our theorem and corollary A. Every concerned L\'evy process shall be written in its canonical form
$$Z_t\; = \; \alpha t \; + \; \int_0^t\int_{|z| \leq 1}\!\!
z\,\mude(ds,dz) \; + \; \int_0^t\int_{|z| > 1}\!\! z\,\mu(ds,dz),$$
and we shall discuss the shape of the jumping measure $\nu$. We refer to Chapter 3 in \cite{sato} for an extensive account on these two examples.

\begin{example}[Stable processes] {\em The measure $\nu$ is given in the integral form
$$\nu(B)\; =\; \int_{\ssc^{d-1}}\lambda(d\xi)\int_0^{+\infty}\Un_B (r\xi)\frac{dr}{r^{1+\beta}}$$
for every measurable set $B\subset\rl^d$, where $0<\beta<2$ and $\lambda$ is some finite positive measure on $\ssc^{d-1}$. We suppose that $\nu$ is non-degenerated, i.e. $\mbox{Supp}\,\lambda$ is not included in any hyper-plane of $\rl^d$. Hence, with the above notations, either $\beta < 1$ and $K =\rl^d$, or $1\leq\beta < 2$ and $L =\rl^d$. It is clear that $Z\in\Wp(0)$ if and only if $p > \beta$. Corollary A reads

\vspace{2mm}

\noindent
(a) If $\beta < 1$, then $Z$ has small deviations in $1$-variation norm if and only if 
$$\alpha\; =\; \frac{1}{1-\beta}\lpa\int_{\ssc^{d-1}}\xi\,\lambda(d\xi)\rpa,$$
i.e. if and only if $Z$ is strictly stable (or the sum of its jumps). 

\vspace{2mm}

\noindent
(b) If $\beta\geq 1$, then $Z$ has small deviations in $p$-variation norm for every $p > \beta$.

\vspace{2mm}

\noindent
Besides, since here $\CC = \CC_{\lambda}$ where $\CC_{\lambda}$ is the convex c\^one generated by $\mbox{Supp}\,\lambda$, one can improve (c) and (d) in Corollary A and show that if $\beta < 1$, then $Z$ has small deviations in $p$-variation norm ($p > 1$) if and only if 
$$\frac{1}{1-\beta}\lpa\int_{\ssc^{d-1}}\xi\,\lambda(d\xi)\rpa\; -\;\alpha\;\in\;\CC_{\lambda}.$$ 

Of course, much more can be said about stable processes. If $\alpha =0$ and $\lambda$ is a symmetric measure (the so-called symmetric $\beta$-stable case) then for every $\gamma > \beta$, the $\gamma$-variation of $Z$ over $[0,1]$ is given by 
$${\lva\lva Z \rva\rva}_{1,\gamma}\; = \; {\lpa\sum_{t\leq 1}{\lva\Delta Z_t\rva}^{\gamma}\rpa}^{1/\gamma}$$
if $\gamma \leq 1$ (notice that this expression makes sense even if $\beta < \gamma < 1$), and satisfies
$${\lva\lva Z \rva\rva}_{1,\gamma}\; \geq \; {\lpa\sum_{t\leq 1}{\lva\Delta Z_t\rva}^{\gamma}\rpa}^{1/\gamma}$$
if $\gamma >1$. Notice that the process
$$S : t\; \longmapsto\; \sum_{s\leq t}{\lva\Delta Z_s\rva}^{\gamma}$$
is a $(\beta/\gamma)$-stable subordinator, whose Laplace transform is given by
$$\esp\lcr\exp -u S_1\rcr\; =\; \exp -\lcr\frac{c_{\lambda}\Gamma(1-\delta)}{\beta}\; u^{\delta}\rcr,$$
where we set $\delta = \beta/\gamma$ and
$$c_{\lambda}\; =\; \int_{\ssc^{d-1}}{\lva\xi\rva}^{\gamma}\lambda(d\xi)$$
(see e.g. Example 24.12. in \cite{sato}). Hence, by De Bruijn's Tauberian theorem (see Theorem 4.12.9 in \cite{bgt}), we get
$$-\log \pb\lcr S_1 < \ee\rcr\; \sim\; (1-\delta){\lpa\frac{c_{\lambda}\Gamma(1-\delta)}{\gamma\ee}\rpa}^{\delta/(1-\delta)}$$
as $\ee\rightarrow 0$. This leads to
$$\lim_{\ee\rightarrow 0}\ee^{\frac{\gamma\beta}{\gamma-\beta}}\log\pb\lcr{\lva\lva Z \rva\rva}_{1,\gamma}\; < \; \ee\rcr\; = \;- \lpa \frac{(\gamma -\beta){\lpa c_{\lambda}\Gamma(1-\beta/\gamma)\rpa}^{\frac{\beta}{\gamma-\beta}}}{\gamma^{\frac{\gamma}{\gamma -\beta}}}\rpa\; =\; - C_{\beta, \gamma}$$
if $\beta < \gamma \leq 1$, and to 
$$\limsup_{\ee\rightarrow 0}\ee^{\frac{\gamma\beta}{\gamma-\beta}}\log\pb\lcr{\lva\lva Z \rva\rva}_{1,\gamma}\; < \; \ee\rcr\;\leq\; -\lpa \frac{(\gamma -\beta){\lpa c_{\lambda}\Gamma(1-\beta/\gamma)\rpa}^{\frac{\beta}{\gamma-\beta}}}{\gamma^{\frac{\gamma}{\gamma -\beta}}}\rpa$$
if $\gamma > 1$.
This can be viewed as small ball probability estimates for symmetric $\beta$-stable processes under the $\gamma$-variation norm, and can probably be extended to more general symmetric L\'evy processes without Gaussian part. In \cite{varsta}, we prove that 
$$\lim_{\ee\rightarrow 0}-\ee^{\frac{\gamma\beta}{\gamma-\beta}}\log\pb\lcr{\lva\lva Z \rva\rva}_{1,\gamma}\; < \; \ee\rcr$$
exists and is finite for {\em every} $\gamma > \beta$. However we could not identify this limit as yet, when $\gamma > 1$.}

\end{example}

\begin{remarks}{\em (a) When $\gamma$ tends to $+\infty$, the constant $C_{\beta,\gamma}$ tends to 1 whereas heuristically, the above limit tends to
$$\lim_{\ee\rightarrow 0}-\ee^{\beta}\log\pb\lcr{\lva\lva Z \rva\rva}_{1,\omega}\; < \; \ee\rcr$$
where ${\lva\lva Z \rva\rva}_{1,\omega}$ stands for the oscillation of $Z$ over $[0,1]$ (see \cite{chisty} Prop. 2.3. p. 27). It follows from the classical result of Taylor \cite{tay1} under the uniform norm (see also \cite{bert} Prop. VIII.3) and standard sublinearity arguments that the latter limit actually exists and belongs to $]0,+\infty[$ (but nothing is known about its explicit value). 

\vspace{2mm}

\noindent
(b) For a linear Brownian motion $W$, it follows from the general result of Stolz \cite{stolz} that 
$$0\; <\;\liminf_{\ee\rightarrow 0}-\ee^{\frac{2p}{p-2}}\log\pb\lcr{\lva\lva W \rva\rva}_{1,p}\; <\; \ee\rcr\;\leq\;\limsup_{\ee\rightarrow 0}-\ee^{\frac{2p}{p-2}}\log\pb\lcr{\lva\lva W \rva\rva}_{1,p}\; <\; \ee\rcr\; <\; +\infty$$ 
for every $p > 2$. This speed of convergence is in accordance with the results of Baldi and Roynette under the H\"older norms \cite{balroy}, and with our previous computation for non-Gaussian symmetric stable processes. It follows from the general results of \cite{lifsim} that the limit actually exists, but we were not able to identify it as yet. As far as we know, the value of the small ball constant for linear Brownian motion is still unknown under the H\"older norms (see the introduction in \cite{balroy}). 

\vspace{2mm}

\noindent
(c) We notice finally that in two celebrated papers \cite{tay2} \cite{fritay}, the {\em exact} variation functions of Brownian motion and stable processes had been computed.}
\end{remarks}

\begin{example}[Self-decomposable processes] {\em The measure $\nu$ is given in the integral form
$$\nu(B)\; =\; \int_{\ssc^{d-1}}\lambda(d\xi)\int_0^{+\infty}\Un_B (r\xi)k_{\xi}(r)\frac{dr}{r}$$
for every measurable set $B\subset\rl^d$, where $\lambda$ is some finite positive measure on $\ssc^{d-1}$, and $k_{\xi}(r)$ is a non-negative function measurable in $\xi$ and decreasing in $r$. 

These class of processes includes the above stable ones, but its range is much wider. In particular it is possible that $K$ and $L$ together are non-trivial: consider for example over $\rl^3 =\{(z_1, z_2, z_3)\}$ endowed with its canonical Euclidean structure,
$$\lambda(d\theta, d\phi)\; =\; \Un_{\{0\leq \theta,\phi\leq\pi/2\}}d\theta\, d\phi\;\;\;\mbox{and}\;\;\; k_{\theta,\phi}(r)\; =\; \frac{\sin\phi\,\cos\theta}{r^{\cos\theta}}\;\Un_{\{r\leq 1\}},$$
where we used the spherical coordinates given by $z_1 = r\sin\phi, z_2 = r \cos\phi\sin\theta$ and $z_3 = r \cos\phi\cos\theta$. Then, with the above notations, we get $K\; =\; \mbox{Vect}\{e_1\},\;L\; =\; \mbox{Vect}\{e_2, e_3\},$ and
$$\Pi^{-1}_K\lpa\AA^{}_K\rpa\; =\; \Pi^{-1}_K\lpa\BB^{}_K\rpa\; =\;\lacc z_1\leq 1 + \pi/2\racc.$$
The associated process $Z$ has finite $p$-variation for every $p > 1$ and since $\CC$ is clearly strictly convex, Corollary A entails that $Z$ has small deviations in $p$-variation norm if 
$$\alpha_1\;\leq\; 1 + \pi/2$$
(notice that here the reverse inclusion is actually true, since $\alpha_1\; >\; 1 + \pi/2$ entails, with the above notation for $Z^1$, that $Z^1_1\; > \; \alpha_1\; -\; \lpa 1 + \pi/2\rpa\;\;\mbox{a.s.}$)

In this example, we stress that even though $\mbox{Dim} \,L \, =\, 2$, there is actually just {\em one} asymptotic direction $u\in L$ as far as our problem is concerned. Indeed, if we set
$$\vtaL\; =\; \Pi_L^{}\lpa\int_{\eta\leq |z|\leq 1}\!\! z\,\nu(dz)\rpa,$$
then we see that $\vtaL = \rta u$ where $u$ is the fixed unit vector  
$$u\; = \; \frac{2e_2 + \pi e_3}{\sqrt{\pi^2 +4}}$$
and where
$$0 \; < \;\rta\;\leq\;\int_{\eta\leq |z|\leq 1}\!\! |z|\,\nu(dz),$$
for every $0 < \eta < 1$. This property, which is not at all a feature of self-decomposable processes (think of $k_{\theta,\phi}(r) = r^{-1} \Un_{\{\theta = 0,\;\phi =0\}} + r^{-3/2} \Un_{\{\theta = 0,\;\phi =\pi/2\}}$ on $\rl^3$), makes the proof of the Theorem somewhat simpler (see Subsection 4.3.1 below).

In any case, we notice that the self-decomposable case is not really relevant to the generality of our result. The monotonicity condition on $k_{\xi}(r)$ entails namely that
$$\CC\; =\; \CC^1\; =\; \CC^{\eta}$$
for every $\eta > 0$, a feature which also simplifies significantly the proof of the Theorem (see the first paragraph in Subsection 4.3.2 below). Actually we have even more: for every $x\in\CC$, there exists $x_1,\ldots, x_d\,\in\,\mbox{Supp}\,\nu$ and $\lambda_1,\ldots, \lambda_d >0$ such that
$$x\; =\; \sum_{i=1}^{d}\lambda_i x_i$$
together with $\mu x_i\in\mbox{Supp}\,\nu$ for {\em every} $\mu > 0$ and $i= 1,\ldots, d$. This latter property makes the proof really easy. The following example, which revisits Example 3, depicts a typical situation where our main result is more difficult to obtain.}
\end{example}

\begin{example}[A pathological measure] {\em Consider on $\rl^+\times\rl^+\times\rl^+$ endowed with the canonical Euclidean structure, the following measure:
$$\nu(dz) \; =\; \Un_{\{0<z_1< z_2^r < z^{rs}_3 < 1\}} z_3^{-2-q}\; dz,$$
where $r,s > 1$ and $s(r+1) + s\; < \; q+1\; < \; s(r +1) + 2\wedge rs$
(notice that $q$ can take any value strictly greater than 2). Then $\nu$ is a jumping measure whose asymptotic subspaces are 
$K\; =\; \mbox{Vect}\{e_1\}$ and $L \; =\; \mbox{Vect}\{e_2, e_3\}$. In this example, the first difficulty  comes from the fact that
$$\frac{\lva\vta_L\rva}{\lva\vta_3\rva}\rightarrow
1,\;\;\;\frac{\lva\vta_3\rva}{\lva\vta_2\rva}\rightarrow
+\infty\;\;\;\mbox{and}\;\;\;\lva\vta_2\rva\rightarrow
+\infty$$
as $\eta\downarrow 0$ (with the obvious notations for $\vta_2$ and $\vta_3$), so that here we must cope with more than one asymptotic direction.  
The second difficulty comes from the degenerescence of $\CC = \lacc z_1 = z_2 = 0,\; z_3\geq 0\racc$. In particular
$$\Pi^{\perp}_3 \lpa\CC\rpa \; =\; \{0\}\;\neq\; \lacc z_1 \geq 0, z_2 \geq 0\racc\; =\; \bigcap_{\eta >0}\overline{\Pi^{\perp}_3 \lpa\Cta\rpa},$$
where $\Pi^{\perp}_3$ stands for the operator of orthogonal projection onto $\mbox{Vect}\{e_1, e_2\}$. This very pathological situation is the matter of Subsection 4.3.2, more particularly of its second paragraph.}

\end{example}

\subsection{Support theorem in $p$-Skorohod topology for Marcus equations}

The principal motivation for our above small deviation result is to prove a support theorem 
\cite{strovar1} for a class of stochastic integral equations driven by $Z$, 
without any assumption on $Z$ but the finiteness of its $p$-variation for some 
$1\leq p < 2$. In this subsection, every index $p$ will be implicitly supposed to belong to $[1,2[$.
We consider on $\rl^m$ 
\begin{eqnarray}
X_t & = & x\; +\;\int_0^t f (X_{s-})\,\diamond\,dZ_s,
\end{eqnarray}
where $x\in\rl^m$ and $f : \rl^m\rightarrow\rl^m\otimes\rl^d$ is a
function which is $\alpha$-Lipschitz for some $p <\alpha < 2$: $f$ is bounded 
with bounded derivatives $\partial_j f$ verifying 
$$\sup_{x\neq y} \frac{\lva\partial_j f(x) - \partial_j f(y)\rva}{{\lva x-y\rva}^{\alpha -1}} \; < \; +\infty.$$
for $j = 1\ldots m$. In (1), the integral is defined followingly:  
$$\int_0^t f\,(X_{s-})\,\diamond\,dZ_s\;=\;\int_0^t f\,(X_{s-})\,dZ_s
\; +\; \sum_{s\leq t}\,g\,(X_{s-},\Delta Z_s),$$ 
where the first integral is a standard It\^o integral and $g :
\rl^m\times\rl^d\rightarrow\rl^m$ is a local Lipschitz function such
that when $(x,z)$ stays in a fixed compact set of $\rl^m\times\rl^d$ 
$$\lva g\,(x,z)\rva\;\leq\; K\,{|z|}^2$$
for some constant $K$, and such that in (1), each time $t$ when $Z$
jumps, $X_t$ is given by the integral in time 1 of the vector field
$x\mapsto f(x)\Delta Z_t$, starting from $X_{t-}$.  

\vspace{2mm}

Introduced by Marcus \cite{marcus}, these stochastic equations are fairly often studied in the literature (see e.g. \cite{fuji} \cite{kpp} \cite{kuni1} \cite{fujikuni} \cite{ishi} \cite{will2} \cite{errami}), even though they concern a specific class of integrand. Their main interest is that they share nice flow properties, and this is not always the case for classical It\^o equations with jumps. Quoting Theorem 7.3.1 in \cite{will1} (which mostly follows from the main result of \cite{lyons1}), we get the following result, which will be the central tool in proving our support theorem: 

\begin{theow}[T.~J.~Lyons, D.~R.~E.~ Williams] Equation (1) can be
  solved path-wise, and has a unique solution. Besides, the map 
$$\Phi : \lacc\begin{array}{l}
           \rl^m\times\Wp\longrightarrow \Wp\\
           (x,Z)\longmapsto X,
           \end{array}
      \right.
$$
is local Lipschitz, where $X$ is the equivalence class of the unique solution to (1) 
starting from $x$, and $\Wp$ is endowed with the $p$-variation norm. 
\end{theow}

\begin{remarks}{\em (a) The local Lipschitz property of $\Phi$ means: for every $T > 0$,  
for every compact set of $\rl^m\times\Wp\lpa [0,T], \rl^d\rpa$ with respect to the norm 
$|.| + {\lva\lva.\rva\rva}_{T,p}$, there exists a constant $K$ such that for every 
$(x,u), (y,v)$ in this compact set,
$${\lva\lva \Phi (x,u) - \Phi(y,v)\rva\rva}_{T,p}\;\leq\; K \lpa |x-y|+{\lva\lva u - 
v\rva\rva}_{T,p}\rpa.$$

\vspace{2mm}

(b) If $Z$ is one-dimensional, then $\Phi : \rl^m\times\sk\longrightarrow \sk$ is actually continuous with respect to the local uniform norm \cite{errami}. Of course, this is no more true when $Z$ is multidimensional and the vector fields defining $f$ do not commute, as it is easily seen by transferring Sussmann's well-known counterexample (see p.~40 in \cite{suss}) to pure jump processes.}

\end{remarks}

We want to find the support of $X$ solution of (1) in $\lpa\sk_p(x), \dd_p\rpa$. Recall that by definition this set is made out of functions $\phi\in\sk_p(x)$ such that for every $n\in\NN^*$ and $\ee > 0$
$$\pb\lcr\, \dd_p^n(X, \phi) < \ee\, \rcr\; >\; 0.$$
As in \cite{marcus} we set $\UU$ for the set of sequences $u = \{u_p\} = \{t_p, z_p\}$, where $\{t_p\}$ is an increasing sequence in $(0,+\infty)$ tending to $+ \infty$ and $\{ z_p\}$ a sequence in Supp $\nu\; -\,\{0\}$. For every $u\in\UU$ and every function $\phi_L : \rl^+\rightarrow L$ with regular $p$-variation, we set 
$$\phi^L_t \; = \;\phi_L(t) \; +\; t\,\alpha_{\nu}$$
for every $t\geq 0$, and introduce the following piecewise differential equation:
\begin{eqnarray}
\psi_t & = & x\; +\;\int_0^t f\,(\psi_s)\,d\phi^L_s\; +\;\sum_{t_p\leq t} 
g_f(\psi^{}_{t_p-}, z_p)
\end{eqnarray} 
where we wrote $g_f(x,z) = f(x)z + g(x,z)$ for every $(x,z)\in\rl^m\times\rl^d$. Notice that since (2) has finitely many jumps on every compact time interval, since $\phi^L$ has regular $p$-variation, and since $f$ is $\alpha$-Lipschitz with $\alpha > p$, the main result of \cite{lyons1} states precisely that there exists a unique solution to (2), which belongs to $\sk_p(x)$. Let $\ssc$ be the set of solutions to (2), $u$ varying in $\UU$ and $\phi_L$ in the set of functions from $\rl^+$ to $L$ with regular $p$-variation. Set $\sso$ for the closure of $\ssc$ in $(\sk_p(x), \dd_p)$. Using together our Theorem and Theorem 7 entails the following 

\begin{corob} Let $f$ be $\alpha$-Lipschitz and $X\in \sk_p(x)$ be the unique solution to $(1)$. Then
$$\mbox{\em Supp}\; X \; = \; \sso.$$
\end{corob}

\begin{remarks}{\em (a) The condition that $f$ is $\alpha$-Lipschitz entails in particular that $f$ is bounded, which is a bit annoying if one wishes to consider e.g. linear equations. In this Corollary we can actually get rid of the boundedness assumption through a standard approximation argument, which we did not include here for the sake of brevity.

\vspace{2mm}

(b) In \cite{marcus} a support theorem was proved in the local Skorohod topology for equation (1) without any assumption on $Z$ and under weaker assumptions on $f$, provided the stochastic part of $Z$ is one-dimensional. H.~Kunita \cite{kuni2} had treated the multidimensional case, but his description of the support is complicated, and his results holds under stringent conditions on the L\'evy measure.

\vspace{2mm}

(c) In both papers \cite{kuni2} and \cite{marcus} the driving process
was allowed to have a Gaussian part. Here we cannot cope with this situation, since Theorem 7 no more holds when the driving process has only finite $p$-variation for some $p\geq 2$. Notice that L\'evy processes without Gaussian part may also have infinite $p$-variation for every $p < 2$ - see Example 2.1. in \cite{will2}. In the Brownian case, Ledoux, Qian and Zhang \cite{lqz} got recently a new proof of Stroock-Varadhan's theorem with the help of Lyons' continuity theorem \cite{lyons2} and a suitable control in $p$-variation ($2<p<3$) on the driving Brownian path together with its L\'evy area process. The same kind of arguments combined with Williams' adaptation of Lyons' theory to jump processes \cite{will1} \cite{will2} could actually be a successful approach to prove the support theorem for Marcus equations in full generality on the L\'evy driving path. However this method promises to be highly technical.

\vspace{2mm}

(d) The unique solution to equation (2) is invariant under $(1+\alpha)$-Lipschitz changes of coordinates - see the final Remarks in \cite{lyons1}. On the other hand, since the integral $\diamond$ is defined through exponentiation of vector fields, (1) is also coordinate-free and may be studied on a nice manifold \cite{fuji} \cite{kuni1}. It is essentially trivial that with an intrinsic definition for the function $g_f$, Corollary B also holds in this more general framework. 

\vspace{2mm}

(e) Viewing $\rl^m$ not as a manifold but as a vector space, classical It\^o equations with jumps are more natural (and more general) objects than Marcus equations. But even in dimension 1 continuity results are not known for such equations, and seem actually quite difficult to prove. We refer however to a recent survey of Dudley and Norvai\v sa \cite{dudnor} for results in this direction in the case of the Dol\'eans-Dade equation. In \cite{ito} a support theorem is obtained in full generality on $Z$ for 1-dimensional It\^o equations, viewing the latter as perturbations of Marcus equations and using a comparison's lemma. This method no more holds in the multidimensional case. See also \cite{supp} for partial results, under heavy assumptions on the L\'evy measure.

\vspace{2mm}

(f) In the literature, there does not seem to exist controllability results in $p$-variation for equation (2). However, considering supports for the local Skorohod topology and using the classical results of \cite{kuni0} we can prove, as in \cite{marcus}, that if $L = \rl^d$
$$\lacc {\rm Lie}_f(y)\; = \; \rl^m\;\;\;\forall\, y\in\rl^m\racc\;\;\Longrightarrow\;\;\lacc\CC_x\;\subset\; \mbox{Supp}\; X\racc$$
where $\CC_x$ stands for the set of continuous functions $\rl^+\rightarrow\rl^m$ starting from $x$ and with the obvious notation for ${\rm Lie}_f(y)$, and that
$$\lacc {\rm Lie}_f(y)\; = \; \rl^m\;\;\forall\, y\in\rl^m\;\;\mbox{and}\;\;\mbox{Supp}\; \nu\; = \; \rl^d\racc\;\;\Longrightarrow\;\;\lacc\mbox{Supp}\; X\; = \; \DD_x\racc$$
where $\DD_x$ stands for the set of c\`ad-l\`ag functions $\rl^+\rightarrow\rl^m$ starting from $x$.
}
\end{remarks}

\section{Some deterministic lemmas}

In this section we gather some easy lemmas about $p$-variation which we will use in proving the Theorem and the Corollaries. We begin with two fairly trivial results:

\begin{lemma} Let $n\in\NN^*$, $T >0$ and $v_0,\ldots, v_n\in\rl^d$. Let $f : [0,T]\rightarrow\rl^d$ be the following step-function: 
$$f(t) = v_i\;\;\;\;\;\mbox{if ${\dsty \;\frac{Ti}{n}\leq t < \frac{T(i+1)}{n}.}$}$$
Then  
$${\lva\lva f\rva\rva}^p_{T, p} \leq n{\max_{0\leq i, j\leq n}\lva v_i - v_j \rva}^p.$$
\end{lemma}
{\em Proof.} Straightforward.
\fin

\begin{lemma} Let $f : \rl^+\rightarrow\rl^d$ be a linear function: $f(t) = ta$ for some $a\in\rl^d$. Then for every $T> 0$, 
$${\lva\lva f\rva\rva}_{T, p} = T\lva a \rva = {\lva\lva f\rva\rva}_{T, \infty}.$$
\end{lemma}
{\em Proof.} Let $0 = t_0 < t_1 < \ldots < t_k = T$ be a partition of $[0,T]$. Writing $s_j = t_j - t_{j-1} > 0$ for $j = 1, \ldots, k$, we get
$$\sum_{j=1}^k{\lva f(t_j) - f(t_{j-1})\rva}^p\; =\;{\lva a \rva}^p \lpa \sum_{j=1}^k s_j^p\rpa\;\leq\; {\lva a \rva}^p{\lpa \sum_{j=1}^k s_j\rpa}^p \; =\;T^p{\lva a \rva}^p$$
since $p\geq 1$. The above inequality is of course an equality when $k = 1$.
\fin

The following definition and lemma will be of constant use in proving the Theorem. The lemma itself is a direct consequence of Lemma 9.

\begin{definition} Let $n\in\NN^*$, $T >0$ and $v\in\rl^d$. The following c\`ad-l\`ag function from $[0,T]$ to $\rl^d$:
$$t\longmapsto \lpa\frac{nt}{T} -k\rpa v\;\;\;\;\;\;\mbox{if $\;{\dsty 
\frac{kT}{n}\leq t < \frac{(k+1)T}{n}}$}$$ 
is called a {\em saw-function} with parameters $(n, T, v)$ and is denoted by $\Saw^{n,T}_v$.
\end{definition}

\begin{lemma} For every $n\in\NN^*$, $T >0$, $v\in\rl^d$ and $p\geq 1$ we have
$${\lva\lva\Saw^{n,T}_v \rva\rva}^p_{T,p}\; =\; 2n{\lva v \rva}^p.$$
\end{lemma}
{\em Proof.} Let $0 = t_0 < t_1 < \ldots < t_k = T$ be a partition of $[0,T]$. 
Writing $q_0 =0$ and 
$$q_j\;=\;\sup\lacc q\in\{0,\ldots, k\}\;/\; t_q < \frac{jT}{n}\racc$$ 
for every $j = 1,\ldots, n+1$, we get
\begin{eqnarray*}
\sum_{i = 1}^k {\lva\Saw^{n,T}_v(t_i) - \Saw^{n,T}_v (t_{i-1})\rva}^p & = & 
\sum_{j = 0}^{n}\sum_{q_j + 1\leq q \leq q_{j + 1}} {\lva\Saw^{n,T}_v(t_q) - 
\Saw^{n,T}_v (t_{q-1})\rva}^p\\
& \leq &\sum_{j = 1}^{n}{\lva\Saw^{n,T}_v(t_{q_j +1}) - \Saw^{n,T}_v 
(t_{q_j})\rva}^p\\
& + & \sum_{j = 0}^{n-1}\sum_{k_j < q \leq q_{j + 1}} {\lva\Saw^{n,T}_v(t_q) - 
\Saw^{n,T}_v (t_{q-1})\rva}^p\\
& \leq &  n{\lva v \rva}^p\; +\; \sum_{j = 0}^{n-1}\sum_{k_j < q \leq q_{j + 1}} 
{\lva\Saw^{n,T}_v(t_q) - \Saw^{n,T}_v (t_{q-1})\rva}^p\\
& \leq &  n{\lva v \rva}^p\; +\;n{\lva v \rva}^p\; = \;2n{\lva v \rva}^p
\end{eqnarray*}
where in the third line we wrote $k_0 = 0$, $k_j = q_j + 1$ for $1\leq j\leq 
n+1$, and used the fact that $q_n + 1 = q_{n+1} = k$, and where in the last 
inequality we used Lemma 9. 

Considering now the following partition of $[0,T]$:
$$0\;<\;\frac{T-\rho}{n}\; <\;\frac{T}{n}\;<\;\frac{2T-\rho}{n}\; 
<\;\frac{2T}{n}\; <\;\cdots\; <\; T - \frac{\rho}{n}\; <\; T,$$
and letting $\rho$ tend to 0, it is easy to see that the above upper bound is 
actually the lowest possible.
\fin

The next elementary lemma, whose statement is trivial if $|a| =0$,  will be used in proving Corollary A.

\begin{lemma} let $\phi$ be a c\`ad-l\`ag function $\rl^+\rightarrow\rl^d$ with finite 1-variations such that for every $t > 0$
$$\phi_t \; =\; t a \; +\; \sum_{s\leq t} \Delta\phi_s,$$
for some fixed vector $a$. Then for every $T >0$
$${\lva\lva\phi\rva\rva}_{T,1}\; =\;T\lva a\rva\; +\; \sum_{t\leq T} \lva\Delta\phi_t\rva.$$
\end{lemma}
{\em Proof.} Fix $T >0$. Introduce 
$$\phi^n : t\; \mapsto\; t a \; +\; \sum_{s\leq t} \Delta\phi_s\Un_{\{\lva\Delta\phi_s\rva \geq 1/n\}},$$
for every $n\geq 1$. Since $\phi^n$ has finitely many discontinuities and reasoning as in Lemma 11, it is clear that
$${\lva\lva\phi^n\rva\rva}_{T,1}\; =\;T\lva a\rva\; +\; \sum_{t\leq T} \lva\Delta\phi^n_t\rva\; \longrightarrow \; T\lva a\rva\; +\; \sum_{t\leq T} \lva\Delta\phi_t\rva$$
as $n\uparrow +\infty$. On the other hand,
$${\lva\lva\phi -\phi^n\rva\rva}_{T,1}\; =\;\sum_{t\leq T} \lva\Delta\phi_t\rva\Un_{\{\lva\Delta\phi_t\rva \leq 1/n\}}\; \longrightarrow \; 0$$
as $n\uparrow +\infty$, which completes the proof of the lemma.

\fin

The next approximation lemma will be used in proving Corollary B. It is probably well-known in the literature on $p$-variation. Nevertheless we give a proof in order to be more complete, even though this is quite tedious.

\begin{lemma} Let $\phi : \rl^+\rightarrow \rl^d$ have regular finite $p$-variation. For every $\ee > 0$, $T > 0$, there exists $n_0\in\NN$ such that for every $n\geq n_0$
$${\lva\lva\phi - \phi^n\rva\rva}_{T,p}\; <\;\ee,$$
where $\phi^n$ is the polygonal approximation of $\phi$ over $[0, T]$ with step 
$T/n$.
\end{lemma}
{\em Proof.} Fix $\ee > 0$ and $T > 0$. Since $\phi$ has regular paths, we can 
find $n_0\in\NN$ such that for every $n\geq n_0$,
$$\sup_{\tiny{\begin{array}{c}
        0=t_0 < \ldots < t_k = T\\
        |t_j - t_{j-1}|\leq T/n
        \end{array}}}\sum_{j=1}^{k}{\lva \phi(t_j) - \phi(t_{j-1})\rva}^p\; < 
\;\frac{\ee}{2^{3p +2}}.$$
Let $n\geq n_0$. We first show that
$$\sup_{\tiny{\begin{array}{c}
        0=t_0 < \ldots < t_k = T\\
        |t_j - t_{j-1}|< T/n
        \end{array}}}\sum_{j=1}^{k}{\lva \phi^n(t_j) - \phi^n(t_{j-1})\rva}^p\; 
< \;\frac{\ee}{2^{2p +1}}.$$
Let $0 = t_0 < t_1 < \ldots < t_k = T$ be a partition of $[0,T]$ such that $|t_j 
- t_{j-1}|< T/n$ for every $j\geq 0$. Write $q^-_0 = q^+_0 = 0$ and set
$$q^-_j\;=\;\sup\lacc q\in\{0,\ldots, k\}\;/\; t_q <s_j \racc, \;\;\; q^+_j\;=\;\inf\lacc q\in\{0,\ldots, k\}\;/\; t_q \geq s_j\wedge T\racc
$$ 
with the notation $s_j = {\dsty \frac{jT}{n}}$, for every $j\geq 0$. Notice that 
$q^+_n = k = q^-_r = q^+_r$ for every $r>n$, and that 
$$s_j \leq t_{q^+_j}\leq t_{q^-_{j+1}} < s_{j+1}$$
for every $j\geq 0$. We get, reasoning as in Lemma 1,
\begin{eqnarray*}
\sum_{i = 1}^k {\lva\phi^n(t_i) - \phi^n(t_{i-1})\rva}^p & \leq & \sum_{j = 
0}^{n-1}{\lva\phi^n (t_{q^-_{j+1}})- \phi^n (t_{q^+_j})\rva}^p
\;+ \;\sum_{j = 1}^{n}{\lva\phi^n (t_{q^+_j}) - \phi^n (t_{q^-_j})  \rva}^p\\
& \leq & \sum_{j = 0}^{n-1}{\lva\phi \lpa s_{j+1} \rpa - \phi \lpa 
s_j\rpa\rva}^p\;+ \;\sum_{j = 1}^{n}{\lva\phi^n (t_{q^+_j}) - \phi^n 
(t_{q^-_j})\rva}^p\\
& \leq & \frac{\ee}{2^{3p +2}}\;+ \; \sum_{j = 1}^{n}{\lva\phi^n (t_{q^+_j}) - \phi^n 
(t_{q^-_j})\rva}^p.
\end{eqnarray*}
Writing $s^T_j = s_j\wedge T$ for every $j\geq 0$ and reasoning again as in Lemma 1, we can control the second term of the right-hand side:
\begin{eqnarray*}
\sum_{j = 1}^{n}{\lva\phi^n (t_{q^+_j}) - \phi^n 
(t_{q^-_j})\rva}^p & \leq & \sum_{j = 1}^{n}2^{p-1}\lpa{\lva\phi^n 
(t_{q^+_j}) - \phi^n (s_j)\rva}^p + {\lva\phi^n \lpa s_j\rpa - \phi^n 
(t_{q^-_j})\rva}^p\rpa\\
& \leq &\sum_{j = 1}^{n}2^{p-1}\lpa{\lva\phi \lpa  
s_{j+1}^T\rpa - \phi\lpa s_j\rpa\rva}^p + {\lva\phi\lpa s_j\rpa - \phi\lpa  
s_{j-1}\rpa\rva}^p\rpa\\
& \leq & \frac{\ee}{2^{2p +2}}\;\;.
\end{eqnarray*}
This yields finally
\begin{eqnarray*}
\sum_{i = 1}^k {\lva\phi^n(t_i) - \phi^n(t_{i-1})\rva}^p & \leq & \frac{\ee}{2^{3p +2}}\lpa 1\;+ \;2^p\rpa\; \leq \; \frac{\ee}{2^{2p +1}}
\end{eqnarray*}
which is the desired result.

\vspace{2mm}

Let now $0 = t_0 < t_1 < \ldots < t_k = T$ be {\em any} partition of $[0,T]$. 
Define $q^-_j$ and $q^+_j$ as above, and
$$J\; = \; \lacc j\geq 0\;\;\mbox{such that}\;\;t_{q^+_{j+1}} > 
t_{q^+_j}\;\;\mbox{or}\;\; j =n\racc.$$
Notice that again
$$s_j \leq t_{q^+_j}\leq t_{q^-_{j+1}} < s_{j+1}$$
if $j\in J$. Writing $\psi^n = \phi - \phi^n$ for simplicity, we get

\vspace{2mm}

${\dsty \sum_{i = 1}^k {\lva\psi^n(t_i) - \psi^n(t_{i-1})\rva}^p \; = \; 
\sum_{j\in J}{\lva\psi^n (t_{q^+_j})- \psi^n (t_{q^-_j})\rva}^p}$

\vspace{-2mm}

\begin{flushright}
${\dsty \; + \; \sum_{j\in J}\lpa\sum_{\tiny{t_{q^+_j}\leq t_q < t_{q+1} \leq 
t_{q^-_{j+1}}}}\!\!\!\!{\lva\psi^n (t_{q+1}) - \psi^n (t_{q})  \rva}^p\rpa.}$
\end{flushright}

On the one hand, for every $j\in J$ and $t_{q^+_j}\leq t_q < t_{q+1} \leq 
t_{q^-_{j+1}}$,
\begin{eqnarray*}
{\lva\psi^n (t_{q+1}) - \psi^n (t_{q})  \rva}^p & \leq &  2^{p-1}\lpa 
{\lva\phi^n (t_{q+1}) - \phi^n (t_{q})  \rva}^p + {\lva\phi(t_{q+1}) - 
\phi(t_{q})  \rva}^p\rpa,\\
\end{eqnarray*}
so that after summation,
\begin{eqnarray*}
 \sum_{j\in J}\lpa\sum_{\tiny{t_{q^+_j}\leq t_q < t_{q+1} \leq 
t_{q^-_{j+1}}}}\!\!\!\!{\lva\psi^n (t_{q+1}) - \psi^n (t_{q})  \rva}^p\rpa & 
\leq & 2^{p-1}\lpa\frac{\ee}{2^{3p +2}}\; + \; \frac{\ee}{2^{2p +1}}\rpa \; \leq 
\; \frac{\ee}{2}\;\; .
\end{eqnarray*}
On the other hand, after summation,
\begin{eqnarray*}
\sum_{j\in J}{\lva\psi^n (t_{q^+_j})- \psi^n (t_{q^-_j})\rva}^p & \leq & 
2^{p-1}\sum_{j\in J}\lpa {\lva\psi^n (t_{q^+_j})- \psi^n (s_j)\rva}^p +  
{\lva\psi^n (s_{j+1}^T)- \psi^n (t_{q^-_{j+1}})\rva}^p\rpa \\
& \leq & 4^{p-1}\sum_{j\in J}\lpa {\lva\phi (t_{q^+_j})- \phi (s_j)\rva}^p\; +\; 
{\lva\phi^n (t_{q^+_j})- \phi^n (s_j)\rva}^p\right. \\
& + & \left.{\lva\phi(s_{j+1}^T)- \phi (t_{q^-_{j+1}})\rva}^p\; +\; 
{\lva\phi^n (s_{j+1}^T)- \phi^n (t_{q^-_{j+1}})\rva}^p\rpa\\
& \leq & 4^{p-1}\lpa\frac{2\ee}{2^{3p +2}}\; + \; \frac{2\ee}{2^{2p +1}}\rpa \; 
\leq \;\frac{\ee}{2}\;\; .
\end{eqnarray*}
Finally, we get
\begin{eqnarray*}
\sum_{i = 1}^k {\lva\psi^n(t_i) - \psi^n(t_{i-1})\rva}^p & \leq & \ee\;\; .
\end{eqnarray*}
\fin

\begin{remark}{\em Francis Hirsch gave me the following counterexample when $\phi$ is continuous and has finite $p$-variation but not regular paths, in the case $p=1$: let $\mu$ be a probability measure over $[0,1]$, singular with respect to Lebesgue measure, and such that  
$$\phi(x)\; =\;\int_0^x\mu(dy)$$
is continuous ($\phi$ is a so-called Lebesgue function). $\phi$ has finite (but not regular finite) $1$-variation. For every $n\geq 1$ we can write
$$\phi^n (x)\; =\;\int_0^x\mu^n(dy)$$
where $\mu^n$ is absolutely continuous with respect to Lebesgue measure. Hence we get
$${\lva\lva\phi - \phi^n\rva\rva}_{1,1}\; =\;{\lva\lva\phi\rva\rva}_{1,1}+ {\lva\lva\phi^n\rva\rva}_{1,1}\; =\; 2.$$}
\end{remark}

\section{Proof of the Theorem}

We first make the general remark that, obviously, it suffices to
consider the situation where the jumps of $Z$ are bounded by $1$, so
that in particular
$$\Zde_t \; =\;\int_0^t\int_{|z|\leq 1}\!\!\!\! z\,\mude(ds, dz)\; 
+\; t\lpa\int_{|z|\leq 1}\!\!\!\! z_K\,\nu(dz)\rpa$$ 
for every $t>0$. We will separate the proof according to ${\rm Dim}\; L$ with
increasing order of difficulty. The arguments are somewhat similar to
those of the Proposition in \cite{dev}, but here the situation is
significantly more complicated because of the $p$-variation norm.

\subsection{${\rm Dim}\; L\; =\; 0$}

This case is obvious since we can take $p = 1$. In particular,
$${\lva\lva\Zde\rva\rva}_{T,1}\; =\; \sum_{t\leq T}{\lva\Delta Z_t\rva},$$
and the Theorem follows readily from the fact that
$$\int_{|z| \leq 1}\!\!\!\! |z|\,\nu(dz)\; <\; \infty.$$

\subsection{${\rm Dim}\; L\; =\; 1$}

Here the situation is more complicated since we must take $p > 1$. We fix $T$ and $\ee > 0$ once and for all. We first set
$$\vta\; =\; \int_{\eta\leq |z|\leq 1}\!\! z\,\nu(dz)\;\;\;\mbox{and}\;\;\;\vtaL\; =\; \int_{\eta\leq |z|\leq 1}\!\! z_L\,\nu(dz).$$
The asymptotic study of $\vtaL$ and its suitable approximation by elements of $\mbox{Supp}\; \nu$ will be actually the central point in the whole proof. 

\vspace{2mm}

Since ${\rm Dim}\; L\; =\; 1$, it clear that for every $\eta$ and $\rho > 0$ there exists $\xtar\in \mbox{Supp}\; \nu$ such that $\lva\xtar\rva < \eta$ and 
$$\sphericalangle\lpa\vtaL,\xtar\rpa\;\;\leq\;\rho,$$
where $\sphericalangle\lpa .\,, .\rpa$ stands for the Euclidean angle between two vectors. Choosing $\eta<\ee/12T$ and $\rho <\eta$ such that $\rho\lva\vtaL\rva <\ee/12T$, we get
$$\lva \vtaL - \atar\xtar\rva\; <\; \ee/6T$$
for some minimizing integer $\atar$. Besides, we can choose a neighborhood $\Vtar$ of $\xtar$, included in $\lacc |z|< \eta\racc$ and small enough, such that
$$\lva \int_{\Vtar}\!\! z\,\nu(dz) - \btar\xtar\rva\; <\; \ee/6T$$
for another minimizing integer $\btar$. Setting
$$\vtar\; =\; \vtaL \; +\; \int_{\Vtar}\!\! z\,\nu(dz)\;\;\;\mbox{and}\;\;\;\gtar\; =\; \atar\; +\; \btar$$
yields
\begin{eqnarray}
\lva \vtar - \gtar\xtar\rva & < & \ee/3T.
\end{eqnarray}
We now introduce the saw-function with parameters $\lpa\gtar,T, -\xtar\rpa$, which we will write $\Saw^{\eta}_{\rho}$ for the sake of simplicity. By Lemma 6,
$${\lva\lva\Saw^{\eta}_{\rho}\rva\rva}^p_{T, p}\; =\; 2\gtar{\lva \xtar \rva}^p.$$
Hence, letting $\eta$ tend to 0,
\begin{eqnarray*}
{\lva\lva\Saw^{\eta}_{\rho} \rva\rva}^p_{T, p} & \sim & 2{\lva \xtar \rva}^{p-1}\lva\vtar\rva\\
& \leq & 2{\lva \xtar \rva}^{p-1}\int_{\Atar}\!\! \lva z\rva\,\nu(dz),
\end{eqnarray*}
where we wrote $\Atar = \lacc z, \; 1 \geq |z|\geq  \lva\xtar\rva/2\racc$. But since  
$$\int_{|z| \leq 1}\!\! {|z|}^p\nu(dz)\; <\; \infty,$$
in the above inequality the right-hand side tends to 0 as $\eta$ tends to 0, and we get 
\begin{eqnarray}
\lim_{\eta\downarrow 0}{\lva\lva\Saw^{\eta}_{\rho}\rva\rva}_{T, p} & = & 0.
\end{eqnarray}
We now come back to the proof of the Theorem. Set 
$$\Btar = \lacc z,\; |z|\leq \eta\racc\cap{\lpa\Vtar\rpa}^c\;\;\;\mbox{and}\;\;\;\Ztard(t)\; =\; \int_0^t\int_{\Btar}\!\!\!\! z\,\mude(ds, dz)\; +\; t\lpa\int_{\Btar}\!\! z_K\,\nu(dz)\rpa$$
for every $t \geq 0$. We obviously have
\begin{eqnarray*}\lim_{\eta\downarrow 0} \lpa\int_{\Btar}\!\! z_K\,\nu(dz)\rpa & = & 0.
\end{eqnarray*}
Hence by Lemma 1, Lemma 9, and the triangle inequality for ${\lva\lva .\rva\rva}_{T, p}$, we see that
\begin{eqnarray}\lim_{\eta\downarrow 0}{\lva\lva\Ztard\rva\rva}_{T, p} & = & 0\;\;\;\mbox{a.s.}
\end{eqnarray}
We next consider the process $\Ztar = \Zde - \Ztard$, which can be written
$$\Ztar(t) \; =\; \int_0^t\int_{{\lpa\Btar\rpa}^c}\!\! z\,\mu(ds, dz)\; -\; t\vtar$$
for every $t >0$. Writing, for each $k = 0\ldots\gtar$, $\sstar(k) = {\dsty \frac{kT}{\gtar}}$, we see that for every $t\in [0,T]$,
$$\Ztar(t) - \Saw^{\eta}_{\rho}(t)\; =\; \pstar(t)\; +\; \ftar(t)$$
where we introduced
$$\pstar(t)\; =\;\sum_{s\leq t} \Delta\Ztar(s) - \sum_{\sstar(k)\leq t} \xtar\;\;\;\;\mbox{and}\;\;\;\;\ftar(t)\; =\;t\lpa\gtar\xtar -\vtar\rpa.$$
By (3) and Lemma 9, we have
$${\lva\lva\ftar\rva\rva}_{T, p}\; < \; \ee/3$$
if $\eta$ was chosen small enough. Introduce ${\lacc
  \Ttar(k),\Utar(k)\racc}_{k\geq 1}$, the successive jumping times and sizes of $\Ztar$. Since $\xtar\in {\mbox {\rm Supp}}\;\nu$ and since $\Vtar\subset{\lpa\Btar\rpa}^c$, we see that the event
$$\lacc\lva\Ttar(k) -\sstar(k)\rva < \lambda, \lva\Utar(k) -\xtar\rva < \lambda,\;\; k= 1,\ldots, \gtar\racc$$
has positive probability for every $\lambda > 0$. But if $\lambda$ is small enough, then on the latter event $\pstar$ is a step function on $[0,T]$ with $2\gtar$ jumps and such that
$$\lva\pstar(t)-\pstar(s)\rva\; <\; (1+\lambda)\lva\xtar\rva$$
for every $s\neq t\in[0,T]$, so that according to Lemma 8
$${\lva\lva\pstar\rva\rva}_{T, p}\; < \; \ee/3$$
if $\eta$ was chosen small enough. Putting everything together leads to
\begin{eqnarray}
\pb\lcr  {\lva\lva\Ztar - \Saw^{\eta}_{\rho} \rva\rva}_{T, p} < 2\ee/3 \rcr & > & 0.
\end{eqnarray}
if $\eta$ was chosen small enough. Using (4), (5), (6), the independence of $\Ztar$ and $\Ztard$ and the triangle inequality for ${\lva\lva .\rva\rva}_{T, p}$, we finally get
$$\pb\lcr  {\lva\lva\Zde\rva\rva}_{T, p} < \ee \rcr\; > \; 0,$$
which finishes the proof in the case ${\rm Dim}\; L\; =\; 1$.

\subsection{${\rm Dim}\; L\; = \; 2$}

We consider this particular situation in order to clarify the exposition - the arguments are analogous in the case ${\rm Dim}\; L\; > \; 2$, but involve heavier notations. The outline of the proof will be roughly the same as in the preceding subsection, except that here the estimate (3) does not hold in general, so that we will need more elements of Supp $\nu$ to approximate $\vtaL$. 

\vspace{2mm}

For each vector $z\in\rl^d$, we will write $z = (x,y)$ according to the unique decomposition $z = x + y$ with $x\in K$ and $y\in L$. Fix an orthonormal basis of $L$. Thanks to the spatial homogeneity of Poisson measures, we first remark that it suffices to consider the case where 
\begin{eqnarray}
\mbox{Supp}\;\nu & \subset & \lacc (x,y),\;\; y_i \geq 0\;\;\forall\, i = 1,2\racc.
\end{eqnarray}
This choice of (strict) positivity will play a crucial r\^ole in the following. Notice, first, that it entails
$$\lva\vtaL\rva\;\rightarrow\; +\infty$$
 as $\eta\downarrow 0$. Choose a subsequence $\{\eta\}$ along which
$$\lim_{\eta\downarrow 0}\frac{\vta}{\lva\vta\rva}\; = \; \lim_{\eta\downarrow 0}\frac{\vtaL}{\lva\vtaL\rva}\; = \; u_1\in{\cal S}^{d-1}.$$
Set $L_1$ for the line generated by $u_1$ and consider the orthogonal
sum $L= L_1\oplus L_2$. Let $\vta_i$ be the projection of $\vtaL$ onto
$L_i$ for $i = 1, 2$. Clearly we have 
$$\lva\vta_1\rva\rightarrow
+\infty\;\;\;\mbox{and}\;\;\;\frac{\lva\vta_1\rva}{\lva\vta_2\rva}\rightarrow
+\infty$$
as $\eta$ tends to 0 along the subsequence. We will consider two
disjoint cases:

\vspace{4mm}

\noindent
{\bf Case A:} There exists a sub-subsequence $\{\eta\}$ along which $\vta_2\rightarrow 0.$

\vspace{2mm}

\noindent
{\bf Case B:} For every sub-subsequence $\{\eta\}$, ${\dsty \liminf_{\eta\downarrow 0}}\lva\vta_2\rva\; > \; 0$ . 

\subsubsection{Case A}

The situation is quite analogous to ${\rm Dim}\; L\; =\; 1$ though a
bit more complicated since here, as we said before, one cannot rely on inequality (3). One
should keep in mind the two-dimensional example where Supp $\nu
\;\subset\;\lacc z = (y_1, y_2),\; y_2 = |y_1|\racc$ and where the
restriction of $\nu$ on each half-line is the same measure. In the
following, each $\eta$ will be implicitly chosen in the sub-subsequence.

\vspace{2mm}

Recall that $\Cta$ stands for the closed convex cone generated by $\Sta\; =\:\mbox{Supp}\;\nu\;\cap\;\lacc z,\; 
|z| <\eta\racc$, and 
$$\CC\; =\;\bigcap_{\eta > 0}\Cta.$$
Notice that clearly $u_1\in\CC$. Besides, since
$$\lim_{\eta\downarrow 0}\frac{\vta_1}{\lva\vta_1\rva}\; =\;
\lim_{\eta\downarrow 0}\frac{\vtaL}{\lva\vtaL\rva}\; =\; u_1,$$
we see that $\vta_1\in\CC$ for $\eta$ small enough. In particular
there exist some integer $r\leq d$ and distinct $\xta_1,\ldots,\xta_r\in\Sta$ such that
\begin{eqnarray}
\lva \vta_1\, - \, \sum_{i=1}^r\ata_i\xta_i\rva\; \leq\; \ee/4T
\end{eqnarray}
for minimizing integers $\ata_1,\ldots,\ata_r$. Notice that by positivity, (7) yields obviously 
\begin{eqnarray}
\lva\ata_i\xta_i\rva & \leq & \lva\vta_1\rva
\end{eqnarray}
for every $i = 1,\ldots, r$. As before we can choose some disjoint neighborhoods $\Vta_i$ of the $\xta_i$, included in $\lacc |z|< \eta\racc$ and small enough, such that
\begin{eqnarray}
\lva \int_{\Vta_i}  z\,\nu(dz) - \bta_i\xta_i\rva & < & \ee/8(r+1)T
\end{eqnarray}
for minimizing integers $\bta_i$. We consider again 
$$\wta_i\; =\; \ata_i\xta_i \; +\; \int_{\Vta_i}\!\! z\,\nu(dz),\;\;\;\;\gta_i\; =\; \ata_i\; +\; \bta_i,$$
and set $\Saw^{\eta}_i$ for the saw-function with parameters
$\lpa\gta_i,T, -\xta_i\rpa$. Using (9) and reasoning exactly as in the
case ${\rm Dim}\; L\; =\; 1$ entail
\begin{eqnarray}
\lim_{\eta\downarrow 0}{\lva\lva\Saw^{\eta}_i\rva\rva}_{T,p} & = & 0
\end{eqnarray}
for every $i = 1, \ldots, r$. Writing
$$\Brta_r = \lacc z,\; |z|\leq
\eta\racc\cap{\lpa\Vta_1\cup\ldots\cup\Vta_r\rpa}^c\;\;\mbox{and}\;\;\Ztad_r(t)\;
=\; \int_0^t\int_{\Brta_r}\!\!\!\! z\,\mude(ds, dz)\; + \; t\lpa\int_{|z|\leq\eta}\!\! z_K\,\nu(dz)\rpa$$
for every $t \geq 0$, we get again, by Lemma 1,
\begin{eqnarray}\lim_{\eta\downarrow 0}{\lva\lva\Ztad_r\rva\rva}_{T, p} & = & 0\;\;\;\mbox{a.s.}
\end{eqnarray}
We finally consider the processes
$$\Zta_i(t) \; =\;\int_0^t\int_{\Vta_i} z\,\mu(ds, dz)\; -\; t\wta_i\;\;\mbox{and}\;\;\Zta_{r+1}(t) \; =\;\int_0^t\int_{\eta\leq |z|\leq 1} z\,\mu(ds, dz)$$
for every $t >0$ and $i= 1,\ldots, r$. Notice that the $\Zta_i$'s are mutually independent and that one can rewrite
$$\Zde(t) \; =\; \Ztad_r (t)\; + \;\sum_{i=1}^{r+1}\Zta_i(t)\; +\; t \lpa\sum_{i=1}^r \ata_i\xta_i  -\vta_1\rpa.$$
for every $t >0$. Using the inequality (10) and reasoning exactly as in the
case ${\rm Dim}\; L\; =\; 1$ yield
\begin{eqnarray}
\pb\lcr  {\lva\lva\Zta_i - \Saw^{\eta}_i \rva\rva}_{T, p} < \ee/4(r+1) \rcr & > & 0
\end{eqnarray}
for every $i = 1, \ldots, r$. Moreover it is clear that 
\begin{eqnarray}
\pb\lcr  {\lva\lva\Zta_{r+1} \rva\rva}_{T, p} < \ee/4(r+1) \rcr & > & 0.
\end{eqnarray}
Using (12), (14), (13), (11), (8), independence arguments and the triangle inequality, we finally get 
$$\pb\lcr  {\lva\lva\Zde\rva\rva}_{T, p} < \ee \rcr\; > \; 0,$$
which completes the proof of the theorem.

\subsubsection{Case B}

This case is the most complicated: here we need to cope with $\vta_2$, a vector which does not belong to $\CC$ in general. One should keep in mind Example 
6. 

\vspace{2mm}

From Case A it is clear that it suffices to prove the following: there exists a fixed integer $r$ and distinct $\xta_1,\ldots,\xta_r\in\Sta$ such that if $\eta$ is chosen small enough along the subsequence
\begin{eqnarray}
\lva \vtaL\, - \, \sum_{i=1}^r\ata_i\xta_i\rva\; \leq\; \ee/2T
\end{eqnarray}
where $\ata_1,\ldots,\ata_r$ are minimizing integers verifying
\begin{eqnarray}
\ata_i{\lva\xta_i\rva}^p & \rightarrow & 0
\end{eqnarray}
for every $i = 1,\ldots, s$, as $\eta$ tends to 0. The estimate (15) would be enough to obtain a small deviation property in uniform norm, as in \cite{dev}. But the latter convergences (16) are crucial to obtain this property in $p$-variation norm, because of Lemma 11. In general (16) is quite difficult to obtain together with (15), since the length of the approximating vectors $\ata_i\xta_i$ is not controlled a priori by that of $\vta_L$, as $\eta\downarrow 0$. Notice, however, that (16) follows readily as soon as   
\begin{eqnarray}
\lva\ata_i\xta_i\rva & \leq & c\lva\vtaL\rva
\end{eqnarray}
for every $i = 1,\ldots, s$ and a constant $c$ independent of $\eta$. The remainder of this article will be devoted to the proof of (15) and (16), a proof which does not require probability theory anymore, but some amount of elementary analysis and geometry.

\vspace{2mm}

Take a sub-subsequence $\{\eta\}$ along which
$$\lim_{\eta\downarrow 0}\frac{\vta_2}{\lva\vta_2\rva}\; = \; u_2\in{\cal S}^{d-1}$$ 
with $\vta_2*u_2 > 0$ for every $\eta$. Set $z_2 = z*u_2$ and $z_2^+ = \sup(0, z_2)$ for every $z\in\rl^d$. Clearly,
\begin{eqnarray}
 \int_{|z|\leq 1} \!\!z_2^+\,\Un_{\{|z_K|< z_2\}}\,\nu(dz) & = & +\infty.
\end{eqnarray}
Let $\Dta$ be the closed convex c\^one generated by $\Sta\cap \{0\, <\, |z_K|\, <\, z_2 \}$ and
$$\DD\; =\;\bigcap_{\eta > 0}\Dta.$$
Set $\Pi_1$ (resp. $\Pi_1^{\perp}$, $\Pi_2$) for the operator of orthogonal projection onto $L_1$ (resp. $L_1^{\perp}$, $L_2$). Because of (18) we see that for every $\eta > 0$,
\begin{eqnarray*}
\frac{{\dsty \Pi_1^{\perp}\lpa \int_{\rho\leq |z|\leq\eta}\!\! z\,\Un_{\{0\, <\, |z_K|\, <\, z_2\}}\,\nu(dz)\rpa}}{{\dsty\lva\Pi_1^{\perp}\lpa\int_{\rho\leq |z|\leq\eta}\!\! z\,\Un_{\{0\, <\, |z_K|\, <\, z_2\}}\,\nu(dz)\rpa\rva}} & \longrightarrow & u_2
\end{eqnarray*}
as $\rho$ tends to 0 along the sub-subsequence. Hence, $u_2\in \overline{\Pi_1^{\perp}\lpa\Dta\rpa}$ for every $\eta > 0$, and 
$$\vta_2\;\in\;\bigcap_{\rho > 0}\overline{\Pi_1^{\perp}\lpa\DD^{\rho}\rpa}$$
for $\eta$ small enough along the sub-subsequence. In particular there exist some integer $s\leq d$ and distinct $\xta_1,\ldots,\xta_s\in\Sta\,\cap\, \{0\, <\, |z_K|\,<\, z_2 \}$ such that
\begin{eqnarray}
\lva \vta_2\, - \, \Pi_1^{\perp}\lpa\sum_{i=1}^{s}\bta_i \xta_i\rpa\rva\; \leq\; \ee/4T
\end{eqnarray}
for minimizing integers $\bta_1,\ldots,\bta_s$. Besides, by positivity, it is clear that
\begin{eqnarray*}
\lva\Pi_2\lpa\bta_i\xta_i\rpa\rva & \leq & \lva\vta_2\rva
\end{eqnarray*}
for every $i = 1,\ldots, s$. Hence, by the triangle inequality,
\begin{eqnarray}
|\Pi_1^{\perp}\lpa\bta_i\xta_i\rpa| & = & \lva\Pi_K\lpa\bta_i\xta_i\rpa\rva + \lva\Pi_2\lpa\bta_i\xta_i\rpa\rva\; \leq\; 2\lva\Pi_2\lpa\bta_i\xta_i\rpa\rva\; \leq\; 2\lva\vta_2\rva
\end{eqnarray}
for every $i = 1,\ldots, s$. We now separate the proof according to the degenerescence of $\DD$ with respect to $L_2$.

\paragraph{The case when $\DD$ is non-degenerated.} We mean the case where 
$$u_2\;\in\;\overline{\Pi_1^{\perp}\lpa\DD\rpa}.$$
Then there exists $c$ independent of $\eta$ such that $u_2 + cu_1\in\DD$, and in (19) $\xta_1,\ldots,\xta_s$ can be chosen such that
\begin{eqnarray}
\lva\Pi_1\lpa\sum_{i=1}^{s}\bta_i\xta_i\rpa\rva & \leq & c\lva\vta_2\rva.
\end{eqnarray}
Since, by positivity,
\begin{eqnarray*}
\lva\Pi_1\lpa\bta_i\xta_i\rpa\rva & \leq & \lva\Pi_1\lpa\sum_{i=1}^{s}\bta_i\xta_i\rpa\rva
\end{eqnarray*}
for every $i = 1,\ldots, s$, we get, by (20) and the triangle inequality,
\begin{eqnarray}
\lva\bta_i\xta_i\rva & \leq & (2+c)\lva\vta_2\rva\;\leq\;(2+c)\lva\vta_L\rva
\end{eqnarray} 
for every $i = 1,\ldots, s$. We now write 
\begin{eqnarray}
\vta_L & = & \lpa\sum_{i=1}^{s}\bta_i\xta_i\rpa\; + \;\lpa\vta_1\, -\, \Pi_1\lpa\sum_{i=1}^{s}\bta_i\xta_i\rpa\rpa\, +\, \lpa \vta_2\, -\, \Pi_1^{\perp}\lpa\sum_{i=1}^{s}\bta_i\xta_i\rpa\rpa
\end{eqnarray}
and set 
$$\wta_1\; = \; \vta_1\, -\, \Pi_1\lpa\sum_{i=1}^{s}\bta_i\xta_i\rpa.$$
Using (21) and the fact that 
$$\lim_{\eta\downarrow 0}\frac{\lva\vta_1\rva}{\lva\vta_2\rva}\; =\; +\infty,$$
we see that there exists $\eta_0$ such that for every $\eta <\eta_0$ along the sub-subsequence,
$$\wta_1*u_1 > 0\;\;\;\mbox{and}\;\;\;\lva\wta_1\rva < \lva\vta_1\rva.$$
Hence, by case A, for every $\eta <\eta_0$ along the sub-subsequence, there exists an integer $r\leq d$ and distinct $\xta_{s+1},\ldots,\xta_{s+r}\in\Sta$ such that 
\begin{eqnarray}
\lva \wta_1\, - \, \sum_{i=1}^{r}\bta_{s+i} \xta_{s+i}\rva\; \leq\; \ee/4T
\end{eqnarray}
where $\bta_{s+1},\ldots,\ata_{s+r}$ are minimizing integers verifying
\begin{eqnarray}
\lva\bta_{s+i}\xta_{s+i}\rva & \leq & \lva\wta_1\rva\; \leq\; \lva\vtaL\rva.
\end{eqnarray}
Clearly together (19), (22), (23), (24) and (25) yield (15) and (17), which completes the proof of the Theorem.

\paragraph{The case when $\DD$ is degenerated.} We mean the delicate situation where 
\begin{eqnarray}
u_2 & \notin & \overline{\Pi_1^{\perp}\lpa\DD\rpa}
\end{eqnarray}
and where in particular (21) no more holds a priori. We denote by $\Dde$ (resp. $\Dtad$) the intersection of $\DD$ (resp. $\Dta$) with $ L = L_1\oplus L_2$. Set $(z_1, z_2)$ for the coordinates on $L_1\oplus L_2$ with respect to $(u_1, u_2)$. Because of (26), $\Dde$ is actually reduced to the half-line $\{z_1 \geq 0, \; z_2 = 0, \; z_K = 0\}$. However, since
\begin{eqnarray*}
u_2 & \in & \bigcap_{\eta > 0}\overline{\Pi_1^{\perp}\lpa\Dta\rpa},
\end{eqnarray*}
for every $\eta >0$ the intersection of $\Dtad$ with the open quadrant $\{z_1 > 0,\; z_2 > 0, \; z_K = 0\}$ is non void. Set $\Delta^{\eta}$ for the frontier of $\Dtad$ in $\{z_1 > 0,\; z_2 > 0, \; z_K = 0\}. \Delta^{\eta}$ is another half-line whose slope with respect to $L_1$ decreases to 0 as $\eta$ decreases to 0. Set $(\zta_1, \zta_2) = (\eta, \zta_2)$ for the point of $\Delta^{\eta}$ with $u_1$-coordinate $\eta$ and consider the increasing convex function 
$$h : \lacc\begin{array}{l}]0,1] \rightarrow \rl^+\\
                           \eta\mapsto\zta_2
           \end{array}
      \right.$$ 
The graph of $h$ is located under a smooth curve with slope $0$ at $\eta = 0$, but we will see that $h$ cannot grow too slowly in the neighborhood of $0$:

\begin{lemma} The function $h :\; ]0,1] \rightarrow \rl^+$ defined above satisfies
$$\int_{|z|\leq 1} \!\! h\lpa |z| \rpa\,\nu(dz)\; =\; +\infty.$$
\end{lemma}

\noindent
{\em Proof.} Since $h$ is positive increasing and since, because of (7), $\mbox{Supp}\,\nu\;\subset\;\lacc z_1 \geq 0\racc$, it suffices to prove that 
$$\int_{|z|\leq 1} \!\! h\lpa z_1\rpa\,\nu(dz)\; =\; +\infty.$$
Suppose first that $d=2$, i.e. $K =\{0\}$. Then clearly, by definition of $h$, 
$$\mbox{Supp}\,\nu\;\cap\;\lacc |z| \leq 1\racc\;\subset\;\lacc z_2^+ \leq h(z_1)\racc.$$
Hence
$$\int_{|z|\leq 1} \!\! h\lpa z_1\rpa\,\nu(dz)\; \geq \; \int_{|z|\leq 1} \!\! z^+_2\,\nu(dz)\; =\;+\infty$$
and this completes the proof of the lemma.

\vspace{2mm}

The case $d > 2$, i.e. $K \neq\{0\}$ is more subtle. First, by the Hahn-Banach theorem, there exists a hyper-plane $H$ containing the half-line $\{z_1 \geq 0, \; z_2 = 0, \; z_K = 0\}$ and separating $\DD$ from $\{z_1 \geq 0, \; z_2 > 0, \; z_K = 0\}$. Its unitary normal vector $n$ oriented in the direction of $\DD$ verifies $n_1 = 0$ and $n_2 < 0$. Besides, we can choose $H$ such that $n_2$ is the lowest possible, in the sense that if $m\in\ssc^{d-1}\cap\lacc z_1 =0\racc$ and if $m_2 < n_2$, then
$$\DD\,\cap\,\lacc m*z < 0\racc\;\neq\;\emptyset.$$
Analogously, for every $\eta > 0$ we set $H^{\eta}$ for the hyper-plane containing $\Delta^{\eta}$, separating $\Dta$ from $\{z_1 > 0,\; z_2 > \atta z_1, \; z_K = 0\}$ (where $\atta = h(\eta)/\eta$), and such that if $n^{\eta}$ is its unitary normal vector oriented in the direction of $\Dta$, then $n^{\eta}_2 < 0$ is the lowest possible. Clearly, we have $n^{\eta}\;\rightarrow n$ and in particular $n^{\eta}_2\;\rightarrow n_2$ as $\eta\downarrow 0$. Hence there exists $\lambda > 0$ and $\eta_0 > 0$ such that $n^{\eta}_2 < -\lambda$ for every $\eta < \eta_0$. A little Euclidean geometry shows then that if $\eta < \eta_0$, then
$$\Dta\,\cap\,\lacc z*n_K \geq 0\racc\;\subset\;\lacc\lambda\lpa\atta z_1 -z_2\rpa + \lva z_K\rva\,\geq\, 0\racc,$$
whereas obviously
$$\Dta\,\cap\,\lacc z*n_K \leq 0\racc\;\subset\;\lacc z_2\,\leq\,\atta z_1 \racc.$$
Hence, for every $\eta < \eta_0$,
\begin{eqnarray*}
\Dta\,\cap\,\lacc \lva z_K\rva \leq \lambda z_2/2\racc & \subset &\lacc z_2\,\leq\, 2\atta z_1 \racc\;\subset\;\lacc z_2\,\leq\, 2 h(z_1) \racc
\end{eqnarray*}
since $\eta\mapsto\atta$ is decreasing. In particular
$$\int_{|z|\leq\eta_0} \!\! z_2\,\Un_{\{0\, <\,|z_K|\,<\, \lambda z_2/2\}}\nu(dz)\; \leq \; 2\int_{|z|\leq\eta_0} \!\!h\lpa z_1\rpa\,\nu(dz)$$ 
But the left-hand side equals $+\infty$ and we get
$$\int_{|z|\leq 1} \!\! h\lpa z_1\rpa\,\nu(dz)\; = \; +\infty,$$
which completes the proof of the lemma.

\fin

\noindent
Set now
$$L(\rho)\; =\;\frac{h(\rho)}{\rho^p}$$
for every $\rho\in\,]0,1]$. Since
$$\int_{|z|\leq 1} \!\! {\lva z\rva}^p\nu(dz)\; < \; +\infty,$$
Lemma 14 obviously entails that
$$\limsup_{\rho\downarrow 0} L(\rho)\; =\; +\infty.$$
In the following we will consider $\lacc\rta\racc$ a sequence in $]0,1]$ with $\rta\leq\eta$, where $\lacc\eta\racc$ is the original sub-subsequence, and such that  
$$L(\rta)\; =\; \sup_{\rta\leq\rho\leq 1} L(\rho)\;\uparrow +\infty.$$
By construction of $h$, we see that for every $\eta > 0$ there exist some integer $s\leq d$, distinct $\xta_1,\ldots,\xta_s\in{\cal S}^{\rta}\,\cap\, \{0\, <\, |z_K|\,<\, z_2 \}$ and minimizing integers $\bta_1,\ldots,\bta_s$ such that
\begin{eqnarray}
\lva \vta_2\, - \, \Pi_1^{\perp}\lpa\sum_{i=1}^{s}\bta_i \xta_i\rpa\rva & \leq & \ee/8T
\end{eqnarray}
and
\begin{eqnarray*}
\rta^{p-1}L\lpa\rta\rpa \frac{{\dsty \lva \Pi_1\lpa\sum_{i=1}^{s}\bta_i \xta_i\rpa\rva}}{{\dsty \lva \Pi_1^{\perp}\lpa\sum_{i=1}^{s}\bta_i \xta_i\rpa\rva}} & \longrightarrow & 1 
\end{eqnarray*}
as $\eta\downarrow 0$. In particular, for every $i = 1,\ldots, s$
\begin{eqnarray*}
\lva\Pi_1\lpa\bta_i \xta_i\rpa\rva &\leq & \lva \Pi_1\lpa\sum_{i=1}^{s}\bta_i \xta_i\rpa\rva\;\;\leq\;\;\frac{2\rta^{1-p}\lva\vta_2\rva}{L\lpa\rta\rpa}
\end{eqnarray*}
and remembering (20), 
\begin{eqnarray*}
\bta_i{\lva\xta_i\rva}^p\;\; \leq\;\; \rta^{p-1}\lpa \lva \Pi_1^{}\lpa\bta_i \xta_i\rpa\rva\, +\, |\Pi_1^{\perp}\lpa\bta_i \xta_i\rpa |\rpa & \leq & 2\lpa \rta^{p-1}\lva\vta_2\rva + \frac{\lva\vta_2\rva}{L\lpa\rta\rpa}\rpa
\end{eqnarray*}
as $\eta\downarrow 0$. On the one hand, since $\rta \leq \eta$,
$$\lim_{\eta\downarrow 0} \rta^{p-1}\lva\vta_2\rva \; =\; 0.$$
On the other hand
\begin{eqnarray*}
\lva\vta_2\rva\;\leq\; \int_{\eta\leq |z|\leq 1}\!\! z_2^+\,\nu(dz) & \leq & \int_{\eta\leq |z|\leq 1}\!\! h(z_1)\,\nu(dz)\;\leq\;\int_{\rta\leq |z|\leq 1}\!\! |z|^p L\lpa |z| \rpa\,\nu(dz).
\end{eqnarray*}
But since $L(\rta)\, =\,{\dsty \sup_{\rta\leq\rho\leq 1}L(\rho)}\,\uparrow\, +\infty$, we have
\begin{eqnarray*}
\int_{\rta\leq |z|\leq 1}\!\! |z|^p\lpa\frac{L\lpa |z| \rpa}{L(\rta)}\rpa\,\nu(dz) & \rightarrow & 0
\end{eqnarray*}
as $\eta\downarrow 0$. This yields
\begin{eqnarray*}
\frac{\lva\vta_2\rva}{L(\rta)} & \rightarrow & 0
\end{eqnarray*}
and, putting everything together,
\begin{eqnarray}
\bta_i{\lva\xta_i\rva}^p & \rightarrow & 0 
\end{eqnarray}
as $\eta\downarrow 0$ along the sub-subsequence, for every $i = 1,\ldots, s$.

\vspace{2mm}

The proof draws now to its final step. Suppose first that
$$\lva \Pi_1\lpa\sum_{i=1}^{s}\bta_i\xta_i\rpa\rva\; \leq \; \lva \vta_1\rva$$
along a subsequence $\{\eta\}$ of the original sub-subsequence. Then if we set
$$\wta_1\; =\;\vta_1\, -\, \Pi_1\lpa\sum_{i=1}^{s}\bta_i\xta_i\rpa$$
as before, this entails that
$$\wta_1*u_1 > 0\;\;\;\mbox{and}\;\;\;\lva\wta_1\rva < \lva\vta_1\rva$$
along this subsequence. Thus, reasoning exactly as above, we can write
$$\vtaL\; =\;\sum_{i=1}^{r+s}\bta_i\xta_i\; +\;\lpa \wta_1\, -\,\sum_{i=1}^{r}\bta_{s+i}\xta_{s+i}\rpa\, +\, \lpa \vta_2\, -\,\Pi^{\perp}_1\lpa\sum_{i=1}^{s}\bta_i\xta_i\rpa \rpa$$
where $r\leq d$ and $\xta_{s+1},\ldots ,\xta_{s+r}$ (resp. $\bta_{s+1},\ldots ,\bta_{s+r}$) are elements of $\Sta$ (resp. minimizing integers) such that (24) and (25) hold. Clearly together (24), (25), (27), and (28) yield (15) and (16), which completes the proof of the Theorem.

\vspace{2mm}

Suppose finally that
\begin{eqnarray}
\lva \Pi_1\lpa\sum_{i=1}^{s}\bta_i\xta_i\rpa\rva & \geq & \lva \vta_1\rva
\end{eqnarray}
along the original subsequence. To treat this very last situation we will remove some mass from $\mbox{Supp}\,\nu$ in the $L_1^+$-direction, digging some subset of $\Sta$, and use Lemma 2 together with Skorohod's absolute continuity theorem - see e.g. Theorem 33.1 in \cite{sato}. Actually, we will need less vectors than in the above situation.

We first appeal to Lemma 2 with $\nu_0 = \nu$ and introduce $\eta_0 > 0$ such that for every $\eta <\eta_0$, every subset $\xita$ of $\{ |z|\leq\eta\}$ and every L\'evy measure $\nub\leq\nu$,
$$\pb\lcr{\lva\lva\Zbtad\rva\rva}_{T, p}\; < \; \ee/2\rcr\; > \; 1/2,$$
where $\mub$ is the Poisson measure over $\rl^+\!\times\rl^d$ with intensity $ds 
\otimes\nub(dz)$, $\mubde = \mub - ds \otimes\nub$, and
$$\Zbtad_t\; = \; \int_0^t\int_{\xita}\!\! z\,\mubde(ds,dz)\; + \; t\int_{|z|\leq\eta}\!\! z_K^{}\,\nu(dz)$$
for every $t > 0$. Choose $\eta < \eta_0$ in the subsequence,  distinct $\xta_1,\ldots,\xta_s\in{\cal S}^{\rta}\,\cap\, \{0\, <\, |z_K|\,<\, z_2 \}$, and minimizing integers $\bta_1,\ldots,\bta_s$ such that (27) and (28) hold. Set
$$\lta\; =\; \inf\{|\xta_1|,\ldots,|\xta_s|\}/2\;\;\;\mbox{and}
\;\;\;\mta\; =\; \lva \Pi_1\lpa\sum_{i=1}^{s}\bta_i\xta_i\rpa\rva\, - \, \lva \vta_1\rva.$$
We may rewrite (27) as 
\begin{eqnarray}
\lva \vtaL\, +\, \mta u_1 \,- \, \sum_{i=1}^{s}\bta_i \xta_i\rva & < & \ee/8T.
\end{eqnarray}
Consider now the restriction of $\nu$ to $\{ |z|\leq\lta\}$. Because of the degeneracy of $\DD$, we see that
for every $\rho > 0$,
\begin{eqnarray}
\int_{|z|\leq\lta}\!\! |z|\,\Un_{\{ {|z|}^2\leq (1 + \rho^2)z_1^2\}}\,\nu(dz) & = & +\infty.
\end{eqnarray}
Take $\rho$ such that $\rho\mta < \ee/8T$. From (31), it is clear that we can find $\lro\in ]0,\lta[$ and a positive measure $\nub\leq\nu$ on $\{ |z|\leq 1\}$, with $\nub$ equivalent to $\nu$ and $\nub$ equal to $\nu$ on $\{|z|\leq \lro\}\;\cup\;\{|z|\geq \lta\}$, such that
\begin{eqnarray}
\lva\uta\; -\; \mta u_1\rva& < & \ee/8T,
\end{eqnarray}
where we set
$$\uta\; =\; \int_{\lro\leq |z|\leq\lta}\!\!\!\!\!\!   z\,\lpa\nu -\nub\rpa (dz)\; =\; \int_{|z|\leq 1}\!\!   z\,\lpa\nu -\nub\rpa (dz).$$
It follows from (30) and (32) that
\begin{eqnarray}
\lva \vtaL\, +\, \uta\,- \, \sum_{i=1}^{s}\bta_i \xta_i\rva & < & \ee/4T.
\end{eqnarray}
Introduce the L\'evy process $\Zba$ given by 
$$\Zba_t\; = \; \int_0^t\int_{|z|\leq 1}\!\! z_K\,\mub(ds,dz)\; +\; \int_0^t\int_{|z|\leq 1}\!\! z_L\,\mubde(ds,dz)\; -\; t\uta$$
for every $t > 0$. Take $\xita = \{ |z|\leq\eta\}\cap{\lpa\Vta_1\cup\ldots\cup\Vta_s\rpa}^c$, where the $\Vta_i$'s are
respective neighborhoods of the $\xta_i$'s in $\{\lta\leq |z|\leq\eta\}$ such that
$$\pb\lcr{\lva\lva\Zbta\rva\rva}_{T, p}\; < \; \ee/2\rcr\; > \; 0,$$
having set
$$\Zbta_t\; = \; \int_0^t\int_{\xita^c}\!\! z\,\mub(ds,dz)\; -\; t\lpa\uta + \int_{\xita^c}\!\! z\,\nub(dz)\rpa$$
for every $t > 0$ (this is clearly possible because of (28), (33), and the reasoning in Case A which led to (13) and (14)). Lemma 2 entails that
$$\pb\lcr{\lva\lva\Zbtad\rva\rva}_{T, p}\; < \; \ee/2\rcr\; > \; 0,$$
so that since $\Zba = \Zbta + \Zbtad$ with $\Zbta$ and $\Zbtad$ independent, we finally get
$$\pb\lcr{\lva\lva\Zba\rva\rva}_{T, p}\; < \; \ee\rcr\; > \; 0$$
by the triangle inequality. But now by Skorohod's absolute continuity theorem, the law of $\Zba$ and $\Zde$ are equivalent for every $\eta > 0$. Hence
$$\pb\lcr{\lva\lva\Zde\rva\rva}_{T, p}\; < \; \ee\rcr\; > \; 0,$$
which completes the proof of the Theorem in the case ${\rm Dim}\; L\; =\; 2$.

\subsection{${\rm Dim}\; L\; >\; 2$}

We briefly describe how this higher dimensional situation can be handled. First, it is clear that we just need to prove (15) and (16) along some subsequence $\{\eta\}$ tending to 0. Again, we can make a choice of strict positivity and suppose that $\mbox{Supp}\, \nu$ is included in a quadrant of $\rl^d$. Take an asymptotic direction $L_1 = \mbox{Vect}\{u_1\}$ and a corresponding subsequence $\{\eta\}$. In order to control the projections of our approximating vectors and to preserve (16), we need to refine our choice of positivity: consider the projection of $\mbox{Supp}\, \nu$ onto $L_1^{\perp}$, take an orthonormal basis of $L_1^{\perp}$ and divide $L_1^{\perp}$ accordingly into $2^{d-1}$ quadrants $Q_1,\ldots, Q_{2^{d-1}}$. Set
$$\nu^{}_{1,i}\; =\; \nu\,\Un_{\{ z_1^{\perp}\in Q_i\}}\;\;\;\mbox{and}\;\;\; \vta_{1,i}\; =\;\int_{\eta\leq |z|\leq 1}\!\! z_1^{\perp}\, \nu^{}_{1,i}(dz)$$
for $i = 1,\ldots, 2^{d-1}$. It is clear that
$$\lim_{\eta\downarrow 0}\frac{\lva\vta_1\rva}{\lva\vta_{1,i}\rva}\; =\; +\infty$$
along the subsequence. Take a sub-subsequence along which either
$$\lim_{\eta\downarrow 0}\frac{\vta_{1,i}}{\lva\vta_{1,i}\rva}\; =\; u^{}_{2,i}\;\in\; L_1^{\perp},\;\;\;\mbox{or}\;\;\;\liminf_{\eta\downarrow 0} \lva\vta_{1,i}\rva\; =\; 0,$$
for every $i = 1,\ldots, 2^{d-1}$, and set $L_{2,i} = \mbox{Vect}\{u^{}_{2,i}\}$. 

Suppose first that ${\rm Dim}\; L\; =\; 3$. Because of our choice of strict positivity for each $\nu^{}_{1,i}$, we can reason as in the situation  ${\rm Dim}\; L\; =\; 2$, Case A or B, and prove that there exists a sub-subsequence $\{\eta\}$, $\xta_{1,i},\ldots,\xta_{r_i,i}\,\in\,\mbox{Supp}\, \nu_{1,i}\,\cap\,\{|z|\,\leq\,\eta\}$, $\bta_{1,i},\ldots,\bta_{r_i,i}$ minimizing integers such that for every $i = 1, \ldots, 2^{d-1}$
\begin{eqnarray*}
\lva \vta_{1,i}\, +\, \rho^{\eta}_{1,i} u^{}_{2,i}\,- \, \Pi_1^{\perp}\lpa \sum_{j=1}^{r_i}\bta_{j,i} \xta_{j,i}\rpa\rva & < & \ee/2^dT
\end{eqnarray*}
where $\rho^{\eta}_{i,1} \geq 0$ and 
\begin{eqnarray*}
\lim_{\eta\downarrow 0} \bta_{j,i} {\lva \xta_{j,i}\rva}^p & = & 0
\end{eqnarray*}
for every $j = 1, \ldots, r_i$. Writing
$$\vtaL\; =\; \vta_1 \; +\; \sum_{i=1}^{2^{d-1}}\vta_{1,i}$$
and reasoning as in the end of Case B leads to
\begin{eqnarray}
\lva \vtaL\, +\, \rta u^{}_1 \, +\, \sum_{i = 1}^{2^{d-1}} \rho^{\eta}_{1,i} u^{}_{2,i}\,- \, \sum_{i=1}^{r}\bta_i \xta_i\rva & < & \ee/2T
\end{eqnarray}
where $\rta\geq 0$ and 
\begin{eqnarray}
\lim_{\eta\downarrow 0} \bta_i {\lva \xta_i\rva}^p & = & 0,
\end{eqnarray}
for some fixed integer $r$, $\xta_i,\ldots,\xta_i\,\in\,\mbox{Supp}\, \nu\,\cap\,\{|z|\,\leq\,\eta\}$, and $\bta_i,\ldots,\bta_r$ minimizing integers. The approximation (34) is not exactly (15), but the (positive) perturbing term
$$\rta u^{}_1 \, +\, \sum_{i = 1}^{2^{d-1}} \rho^{\eta}_{1,i} u^{}_{2,i}$$
can be canceled as above by an absolutely continuous transformation of the law of $\Zde$, after removing some mass from $\mbox{Supp}\, \nu_{1,i}\,\cap\,\{|z|\,\leq\,\lta\}$ for each $i = 1, \ldots, 2^{d-1}$, with $\lta = \inf \{ |\xta_i|\}/2$. This transformation leads to the small deviation property for the original process $\Zde$.

When  ${\rm Dim}\; L$ gets higher, we need to refine again and again our decomposition of $\mbox{Supp}\, \nu$, dividing first each orthogonal of $u^{}_{2,i}$ into $2^{d-2}$ quadrants $R_1,\ldots, R_ {2^{d-2}}$ and introducing
$$\nu^{}_{2,i,j}\; =\; \nu\,\Un_{\{ z_1^{\perp}\in Q_i, z_{2,i}^{\perp} \in R_j\}}\;\;\;\mbox{and}\;\;\; \vta_{2,i, j}\; =\;\int_{\eta\leq |z|\leq 1}\!\! z_{2,i}^{\perp}\, \nu^{}_{2,i,j}(dz)$$
for $i = 1,\ldots, 2^{d-1}$, $j = 1,\ldots, 2^{d-2}$...etc. Setting $k = {\rm Dim}\; L$ and writing
$$\vtaL\; =\; \vta_1 \; +\; \sum_{i=1}^{2^{d-1}}\vta_{2,i}\; +\; \ldots\; +\; \lpa\sum_{i_1 =1}^{2^{d-1}}\ldots\sum_{i_{k-2}=1}^{2^{d-(k-2)}}\vta_{k-1,i_1, \ldots, i_{k-2}}\rpa$$
leads to an approximation of type (34) together with the control (35), which finishes the proof of the Theorem.

\section{Proof of the Corollaries}

\subsection{Small deviations around continuous curves}

The following proposition shows that the small deviation property for $\Zde$ also holds around $L$-valued curves with finite regular $p$-variation. Of course, this would be a direct consequence of the Theorem if one had some kind of Cameron-Martin formula for $\Zde$ as for Brownian motion. But here $\Zde$ has no Gaussian part and it is well-known, for example, that the law of $\Zde^u : t \mapsto\Zde_t + tu$ is not absolutely continuous (and even singular) with respect to the law of $\Zde$ if $u\neq 0$ - see again Theorem 33.1. in \cite{sato}. To prove this proposition we will need Lemma 13 as well as a slight perturbation of the Poisson measure, which replaces in some sense the density transformation.  

\begin{proposition} Let $1 \leq p < 2$ and $Z$ be a L\'evy process with finite 
$p$-variation and parameters $(\alpha, \nu)$. For every $\ee > 0$, $T > 0$ 
and $\phi_L: \rl^+ \rightarrow L$ with finite regular $p$-variation over compact 
sets,
$$\pb\lcr{\lva\lva\Zde - \phi_L\rva\rva}_{T, p}\; < \; \ee\rcr\; > \; 0.$$
\end{proposition}
{\em Proof.} Clearly, we can suppose that $\phi_L(0) =0$ and that the
jumps of $Z$ are bounded by 1. In particular
$$\Zde_t \; =\; \int_0^t\int_{|z|\leq 1}\!\!\!\! z\,\mude(ds, dz)\; +\; t\int_{|z|\leq 1}\!\!\!\! z_K\,\nu(dz)$$
for every $t\geq 0$. Fix 
$\ee > 0$, $T > 0$ and $\phi_L: \rl^+ \rightarrow L$ with finite regular 
$p$-variation over compact sets. By Lemma 13, there exists $n_0\in\NN$ such that 
for every $n\geq n_0$
$${\lva\lva\phi_L - \phi^n_L\rva\rva}_{T, p}\; <\;\frac{\ee}{3},$$
where $\phi^n_L$ is the polygonal approximation of $\phi_L$ over $[0, T]$ with 
step $T/n$. Fix $n\geq n_0$. Let $v_0 = 0$ and $v_1,\ldots, v_n$ be the vectors 
of $L$ defining $\phi^n_L$:
$$\phi^n_L (t)\; =\; \frac{T}{n}(v_0 + \ldots + v_j) + 
(t-s_j)v_{j+1}\;\;\;\;\;\mbox{if $s_j\leq t\leq s_{j+1}$}$$
where again we set ${\dsty s_j = \frac{jT}{n}}$. Set, for every $j = 0,\ldots, 
n$ and $s_j\leq t\leq s_{j+1}$,
$$\Zjd_t \; =\; \Zde_t \; -\; \Zde_{s_j}\; -\; (t- s_j)v_{j+1}$$
and
\begin{eqnarray*}
\Yjd_t & = &\int_{s_j}^t\int_{|z|\leq 1}\!\!\!\! z\,\mude_j(ds,dz)\; +\; 
(t-s_j)\int_{|z|\leq 1}\!\!\!\! z_K\,\nu_j(dz),
\end{eqnarray*}
where we wrote  
$$\nu_j(dz)\; =\;\nu (dz) + 2\lva v_{j+1}\rva\delta_{\!\!\frac{v_{j+1}}{2\lva v_{j+1}\rva}}(dz)$$ 
for every  $j = 0,\ldots, n$ (with the notation $\frac{v_j}{|v_j|} = 0$ if 
$|v_j| = 0$), and where $\mude_j$ is the compensated measure of $\mu_j$, the Poisson measure with
intensity $ds\otimes\nu_j$. Notice that clearly, 
\begin{eqnarray*}
\Yjd_t & = &\int_{s_j}^t\int_{|z|\leq 1}\!\!\!\! z\,\mude_j(ds,dz)\; +\; 
(t-s_j)\int_{|z|\leq 1}\!\!\!\! z_K\,\nu(dz),
\end{eqnarray*}
so that for every $0 < \eta < 1/2$ and $j =0,\ldots, n$,
$$\lacc{\lva\lva\Yjd\rva\rva}_{[s_j,s_{j+1}], p}\; < \; \eta\racc\;\subset\; \lacc \Yjd_t = \Zjd_t\;\;\forall t\in [s_j,s_{j+1}]\racc.$$ 
In particular
$$\lacc{\lva\lva\Yjd\rva\rva}_{[s_j,s_{j+1}], p}\; < \; \eta\racc\; 
=\;\lacc{\lva\lva\Zjd\rva\rva}_{[s_j,s_{j+1}], p}\; < \; \eta\racc.$$
Now, by the Theorem, 
$$\pb\lcr{\lva\lva\Yjd\rva\rva}_{[s_j,s_{j+1}], p}\; < \; \eta\rcr\; > \; 0$$
for every $\eta > 0$ and $j =0,\ldots, n$. Hence we get, by independence of the increments of $\Zde$, 
\begin{eqnarray}
\pb\lcr{\lva\lva\Zjd\rva\rva}_{[s_j,s_{j+1}], p} \; < \; \eta\;\;\;\mbox{for 
every $j =0,\ldots, n$}\;\rcr & > & 0
\end{eqnarray}
for every $0 < \eta < 1/2$. Introduce now the following function:
$$\tphi^n_L (t)\; =\;\phi^n_L (t)\; +\;\sum_{k=0}^{j-1} \Zde^k_{s_{k+1}}\;\;\;\;\;\mbox{if $s_j\leq t < s_{j+1}$},$$
for every $j =0,\ldots, n$. $\tphi^n_L$ is a discontinuous
perturbation of $\phi^n_L$ such that $\tphi^n_L -\phi^n_L$ is a
step-function and $\tphi^n_L(s_j) = \Zde_{s_j}$ for 
every  $j =0,\ldots, n$. On the one hand, reasoning as in Lemma 13, we can choose $\eta$ sufficiently small such that 
\begin{eqnarray}
\lacc{\lva\lva\Zjd\rva\rva}_{[s_j,s_{j+1}], p}\; < \; \eta\;\;\;\mbox{for 
every $j =0,\ldots, n$}\;\racc &\subset &\lacc{\lva\lva\Zde - 
\tphi^n_L\rva\rva}_{T, p}\; < \;\frac{\ee}{3}\racc.
\end{eqnarray}
On the other hand, since according to Lemma 8
$${\lva\lva\phi^n_L - \tphi^n_L\rva\rva}^p_{T, p}\;\leq\; n^{p+1}{\max_{0\leq k\leq n-1}\lva\Zde^k_{s_{k+1}}\rva}^p,$$
we also have
\begin{eqnarray}
\lacc{\lva\lva\Zjd\rva\rva}_{[s_j,s_{j+1}], p}\; < \; \eta\;\;\;\mbox{for 
every $j =0,\ldots, n$}\;\racc &\subset &\lacc{\lva\lva\phi^n_L - 
\tphi^n_L\rva\rva}_{T, p} \; < \;\frac{\ee}{3}\racc
\end{eqnarray}
for $\eta$ small enough. Putting (36), (37), (38) together and using the triangle inequality  complete the proof of the
Proposition.
 
\fin

\subsection{Proof of Corollary A}

(a) By the Theorem we just need to prove the reverse inclusion. Suppose $\alpha_{\nu}\neq 0$. Since
$$Z_t\; =\;\alpha_{\nu}\; +\;\sum_{s\leq t}\Delta Z_s$$
for every $t \geq 0$, we see by Lemma 12 that
\begin{eqnarray*}
\pb\lcr{\lva\lva Z\rva\rva}_{1,1}\; < \; \ee\rcr & = & 0
\end{eqnarray*}
as soon as $\ee < \lva\alpha_{\nu}\rva$.

\vspace{2mm}

\noindent
(b) This follows readily from Proposition 15.

\vspace{2mm}

\noindent
(c) $\,$ Fix $\ee, T > 0$. Since $\alpha\in \Pi^{-1}_K\lpa\AA^{}_K\rpa$, there exists $\alpha^{}_L\in L$ such that
$$\beta\; =\; \alpha^{}_L -\alpha_{\nu}\;\in\;\CC.$$
Hence, for every $\eta >0$, there exists $\xta_1,\ldots,\xta_d\in\mbox{Supp}\,\nu\;\cap\;\{|z|\leq\eta\}$ and $\ata_1,\ldots,\ata_d$ 
minimizing integers, such that
\begin{eqnarray}
\lva \beta\, - \, \sum_{i=1}^d\ata_i\xta_i\rva\; \leq\; \ee/4T.
\end{eqnarray}
Besides, since $\CC$ is strictly convex, it is clear that there exists $c > 0$ independent of $\eta$ such that
\begin{eqnarray}
\lva\ata_i\xta_i\rva & \leq & c\lva\beta\rva
\end{eqnarray}
for every $i =1,\ldots,d$. Introduce now $\rta\; =\;\inf\{|\xta_i|,\; i = 1,\ldots, d\}/2$, and decompose $Z$ into
$$Z\; =\; \Ztad\; +\Zta$$
where we set
$$\Ztad_t\; = \;\int_0^t\int_{|z| \leq \rta}\!\! z^{}_K\,\mu(ds,dz)\; +\;\int_0^t\int_{|z| \leq \rta}\!\! z^{}_L\,\mude(ds,dz)\; +\; t\lpa\alpha^{}_L \, -\,\int_{\rta\leq |z| \leq 1}\!\! z^{}_L\,\nu(dz)\rpa$$
and
$$\Ztad_t\; = \;\sum_{s\leq t}\Delta Z_s \Un_{\{\lva\Delta Z_s\rva >\rta\}}\; -\; t\beta$$
for every $t > 0$. The Theorem and Proposition 15 yield readily
$$\pb\lcr{\lva\lva\Ztad\rva\rva}_{T,p}\; < \; \ee/2\rcr\; > \; 0.$$
Hence, by independence and the triangle inequality, it suffices to show that
$$\pb\lcr{\lva\lva\Zta\rva\rva}_{T,p}\; < \; \ee/2\rcr\; > \; 0.$$
It is now clear that the latter can be done through (39), (40), and the same approximation procedure which we used repeatedly during the proof of the Theorem. 

\vspace{2mm}

\noindent
(d) This follows readily from the main Theorem in \cite{dev} and from the inequality
$$\pb\lcr{\lva\lva\Zta\rva\rva}_{T,p}\; < \; \ee \rcr\; \leq \;\pb\lcr{\lva\lva\Zta\rva\rva}_{T,\infty}\; < \; \ee\rcr$$
for every $T, \ee >0$.

\subsection{Proof of Corollary B}

We first quote a lemma which is a direct consequence of Lyons' continuity theorem \cite{lyons1}.

\begin{lemma} Let ${\lacc x^i_t, \;0\leq t\leq T\racc}_{i=1,2}$ be the solutions to the following rough differential equations on $\rl^m$:
$$ x^i_t \;  = \;x_i\; +\;\int_0^t \!\! f\lpa x^i_s\rpa\,dz_s,$$  
where $z$ is a function with regular finite $p$-variation and $f$ an
$\alpha$-Lipschitz vector field with $\alpha > p$. Then there exists a constant $K$ (depending on $T$ and $f$) such that
$${\lva\lva x^1 - x^2\rva\rva}_{T, p}\; \leq\; K\lva x_1 -x_2\rva.$$
\end{lemma}
We can now proceed to the proof of Corollary B, which will mimic that of the Theorem in \cite{marcus}. 
The first inclusion $\mbox{Supp}\; X \; \subset \; \sso$ is an easy consequence
of the fact that for every $n\geq 1$
$$\lim_{\eta\rightarrow 0}{\lva\lva X - \Xta \rva\rva}_{n, p}\; =\; 0$$
where $\Xta$ is the solution to (1) with $\nu$ replaced by $\Un_{|z|\geq\eta}\nu(dz)$ - which follows readily from Lemma 1 and Theorem 7, and of the usual routine which may be found e.g. in \cite{strovar1}.

\vspace{2mm}

The second inclusion $\sso\;\subset\;\mbox{Supp}\; X$ will be a consequence of Theorem 7 and Proposition 15, as in \cite{marcus}. Fix $n\in\NN^*$, $\ee > 0$, $u\in\UU$, and $\phi_L : \rl^+\rightarrow L$ with regular $p$-variation. Let $\psi$ be the solution of (2) given by $u$ and $\phi_L$. Let $N_n\in\NN^*$ be such that
$$t_0 = 0 < t_1 < \ldots < t_{N_n} \leq n + 1 < t_{N_n + 1}$$
are the successive jumping times of $\psi$. Introduce
$$\eta\; =\;\inf\{|z_i|,\; i = 1,\ldots, N_n\}/2\;\;\;\mbox{and}
\;\;\;\Zta_t\; = \;\int_0^t\int_{|z|\geq\eta}\!\! z\; \mu(ds,dz)$$
for every $t\geq 0$. Set $\{T_q\}$ for the sequence of $\Zta$'s successive
jumping times, and $\tpsi$ for the solution of (2) where
$\{t_q\}$ is replaced by $\{T_q\}$. For every $\rho > 0$, the event
$$\lacc\sup_{1\leq q\leq N_n +1}\lva T_q -t_q\rva < \rho\racc$$
has a positive probability. We now introduce $\lambda$, the only piecewise linear change of time transforming $t_q$ into $T_q$ for each $q=1, \ldots, N_n$, and whose right derivative takes its values in $\lacc 1/2, 1, 2\racc$. Thanks to the continuity of $\psi$ (resp. of $\tpsi$) on each $]t_i, t_{i+1}[$ (resp. on each $]T_i, T_{i+1}[$) and to a repeated use of Lemma 16, it is easy to see that
$$\lacc\sup_{1\leq q\leq N_n +1}\lva T_q -t_q\rva < \rho\racc\; \subset\; \lacc{\lva\lva\psi\circ\lambda-\tpsi\rva\rva}_{[0,n+1], p}\; < \;\ee/2\racc$$
for $\rho > 0$ small enough. Hence
$$\lacc\sup_{1\leq q\leq N_n +1}\lva T_q -t_q\rva < \rho\racc\; \subset\; \lacc\dd^n_p\lpa\psi,\tpsi\rpa\; < \;\ee/2\racc$$
for $\rho > 0$ small enough. On the other hand, Proposition 15 entails easily that
$$\pb\lcr{\lva\lva\Ztad - \phi^L\rva\rva}_{n+1, p}\; < \; \rho\rcr\; > \; 0$$
for every $\rho > 0$, where we set $\Ztad = Z - \Zta$. Using now Theorem 7 and reasoning exactly as in the proof of the Theorem of \cite{marcus} (under the $p$-variation norm) entail
$$\pb\lcr{\lva\lva X - \tpsi\rva\rva}_{n+1, p}\; < \; \ee/2,\;\;\dd^n_p\lpa\psi,\tpsi\rpa\; < \;\ee/2\rcr\; > \; 0,$$
which finishes the proof since obviously
$$\lacc{\lva\lva X - \tpsi\rva\rva}_{n+1, p}\; < \; \ee/2,\;\;\dd^n_p\lpa\psi,\tpsi\rpa\; < \;\ee/2\racc\; \subset \; \lacc\dd^n_p\lpa X,\psi\rpa\; < \;\ee\racc.$$

\vspace{2mm}

\noindent
{\em Acknowledgements:} I thank Yasushi Ishikawa, Michel Lifshits and an anonymous referee for several fruitful comments on the first drafts of this paper.

\vspace{10mm}

\normalsize
\noindent
\'Equipe d'Analyse et Probabilit\'es\\
 Universit\'e d'\'Evry-Val 
d'Essonne\\
Boulevard Fran\c{c}ois Mitterrand\\ 
F-91025 \'EVRY cedex\\
e-mail: {\tt simon@maths.univ-evry.fr}  


\footnotesize

\begin{thebibliography}{10}

\vspace{2mm}

\bibitem{balroy}
P.~Baldi and B.~Roynette.
\newblock
Some exact equivalents for the Brownian motion in H\"older norm. 
\newblock {\em Probab. Th. Rel. Fields} {\bf 93} (4), pp. 457-484, 1992.

\bibitem{bgt}
N.~H.~Bingham, C.~M.~Goldie and J.~L.~Teugels.
\newblock {\em Regular Variation}.
\newblock Cambridge University Press, Cambridge, 1987.

\bibitem{bert}
J.~Bertoin.
\newblock {\em L\'evy processes}.
\newblock Cambridge University Press, Cambridge, 1996.

\bibitem{breta}
J.~Bretagnolle.
\newblock $p$-variation de fonctions al\'eatoires I et II.
\newblock In: {\em S\'em. Prob.} {\bf VI}, pp. 51-71. Springer, Berlin, 1971.

\bibitem{chisty}
V.~V.~Chistyakov and O.~E.~Galkin.
\newblock On maps of bounded $p$-variation with $p > 1$.
\newblock {\em Positivity} {\bf 2} (1), pp. 19-45, 1998.

\bibitem{dudnor}
R.~M.~Dudley and R.~Norvai\v sa.
\newblock 
An introduction to $p$-variation and Young integrals. With emphasis on sample functions of stochastic processes. 
\newblock {\em MaPhySto Lect. Notes} {\bf 1}, Univ. of Aarhus, 1998.

\bibitem{errami}
M.~Errami, F.~Russo and P.~Vallois.
\newblock It\^o formula for ${\CC^{1,\lambda}}$ functions of a c\`adl\`ag
  process and related calculus 
\newblock {\em Prob. Theory and Rel. Fields} {\bf 122} (2), pp. 191-221, 2002.

\bibitem{fritay}
B.~E.~Fristedt and S.~J.~Taylor.
\newblock Strong variation for the sample function of a stable process. 
\newblock {\em Duke Math. J.} {\bf 40} (1), pp. 259-278, 1973.

\bibitem{fuji}
T.~Fujiwara.
\newblock Stochastic differential equations of jump type on manifolds and 
L\'evy flows. 
\newblock {\em J. Math. Kyoto Univ.} {\bf 31} (1), pp. 99-119, 1991.

\bibitem{fujikuni}
T.~Fujiwara and H.~Kunita.
\newblock Canonical stochastic differential equations based on
semimartingales with spatial parameters I and II.
\newblock {\em Kyushu J. Math.} {\bf 53}, pp. 265-331, 1999.

\bibitem{ishi}
Y.~Ishikawa.
\newblock Existence of the density for a singular jump process and its short time properties.
\newblock {\em Kyushu J. Math.} {\bf 55}, pp. 267-299, 2001.

\bibitem{jacshy}
J.~Jacod and A.~N. Shiryaev.
\newblock {\em Limit theorems for stochastic processes}.
\newblock Springer-Verlag, Berlin, 1987.

\bibitem{kuni0}
H.~Kunita.
\newblock Supports of diffusion processes and controllability problems.
\newblock In: {\em Proceedings of the International Symposium on Stochastic 
Differential Equations (Kyoto 1976)}, pp. 163-185. Wiley, Brisbane, 1978.

\bibitem{kuni1}
H.~Kunita.
\newblock Stochastic differential equations with jumps and stochastic
flows of diffeomorphisms. 
\newblock In: Ikeda et al. (editor), {\em It\^o's stochastic calculus
and probability theory}, pp. 197-211. Springer, Berlin, 1996.

\bibitem{kuni2}
H.~Kunita.
\newblock Canonical stochastic differential equations based on {L}\'evy
  processes and their supports.
\newblock In Crauel et al. (editor), {\em Stochastic dynamics},
pp. 283-304. Springer, New York, 1999.

\bibitem{kpp}
T.~G. Kurtz, E.~Pardoux and Ph.~Protter.
\newblock Stratonovich stochastic differential equations driven by general
  semimartingales.
\newblock {\em Ann. I.H.P. Prob. Stat.} {\bf 31}, pp. 351-377, 1995.

\bibitem{lqz}
M.~Ledoux, Z.~Qian and T.~Zhang.
\newblock Large deviations and support theorem for diffusion processes via rough paths.
\newblock {\em Stoch. Proc. Appl.} {\bf 102} (2), pp.265-283, 2002.

\bibitem{lifsim}
M.~A.~Lifshits and T.~Simon.
\newblock Small deviations for fractional stable processes.
\newblock Submitted to {\em Annales de l'Institut Henri Poincar\'e}, 2003.

\bibitem{lyons1}
T.~J.~Lyons.
\newblock Differential equations driven by rough signals (I): An extension of an 
inequality of L.~C.~Young.
\newblock {\em Math. Res. Letters} {\bf 1} (4), pp. 451-464, 1994.

\bibitem{lyons2}
T.~J.~Lyons.
\newblock Differential equations driven by rough signals.
\newblock {\em Rev. Mat. Iberoam.} {\bf 14} (2), pp. 215-310, 1999.

\bibitem{marcus1}
S.~I.~Marcus.
\newblock Modelling and approximation of stochastic differential equations 
driven by semi-martingales. {\em Stochastics} {\bf 4}, pp. 223-245, 1981.

\bibitem{sato}
K.-I.~Sato.
\newblock {\em L\'evy processes and infinitely divisible distributions.}
\newblock Cambridge University Press, Cambridge, 1999.

\bibitem{sharpe}
M.~J.~Sharpe.
\newblock Support of convolution semi-groups and densities.
\newblock In: Heyer (editor), {\em Probability measures on groups and
related structures {\bf XI}}, pp. 364-369. World Scientific
Publishing, River Edge, 1995.

\bibitem{supp}
T.~Simon.
\newblock Support theorem for jump processes.
\newblock {\em Stoch. Proc. Appl.} {\bf 89} (1), pp. 1-30, 2000.

\bibitem{marcus}
T.~Simon.
\newblock Support of a Marcus equation in dimension 1.
\newblock {\em Elec. Comm. Probab.} {\bf 5}, pp. 149-157, 2000.

\bibitem{dev}
T.~Simon.
\newblock Sur les petites d\'eviations d'un processus de
{L}\'evy. {\em Pot. Analysis} {\bf 14} (2), pp. 155-173, 2001. 

\bibitem{ito}
T.~Simon.
\newblock Support d'une \'equation d'It\^o en dimension 1. In: {\em S\'em. Prob.} {\bf 36}, pp. 314-330, 2002. 

\bibitem{varsta}
T.~Simon.
\newblock Small ball estimates in $p$-variation for stable processes. In preparation, 2003.

\bibitem{stolz}
W.~Stolz.
\newblock Une m\'ethode \'el\'ementaire pour l'\'evaluation de petites boules browniennes. {\em C. R. Acad. Sci. Paris S\'er. I} {\bf 316} (11), pp. 1217-1220, 1993. 

\bibitem{strovar1}
D.~W. Stroock and S.~R.~S. Varadhan.
\newblock On the support of diffusion processes with applications to the strong
  maximum principle.
\newblock In: {\em Proc. 6th Berkeley Symp. Math. Stat. Prob.} {\bf III},
  pp. 333-359. Univ. California Press, Berkeley, 1972.

\bibitem{suss}
H.~J.~Sussmann.
\newblock On the gap between deterministic and stochastic ordinary differential equations.
\newblock {\em Ann. Prob.} {\bf 6}, pp. 19-41, 1978.

\bibitem{tay1}
S.~J.~Taylor.
\newblock Sample path properties of a transient stable process.
\newblock {\em J. Math. Mech.} {\bf 16}, pp. 1229-1246, 1967.

\bibitem{tay2}
S.~J.~Taylor.
\newblock Exact asymptotic estimates of Brownian path variation. 
\newblock {\em Duke Math. J.} {\bf 39} (1), pp. 219-241, 1972.

\bibitem{tortr1} 
A.~Tortrat. 
\newblock Le support des lois ind\'efiniment divisibles dans un groupe 
Ab\'elien localement compact. 
\newblock {\em Math. Z.} {\bf 197}, pp. 231-250, 1988.

\bibitem{will1}
D.~R.~E.~Williams.
\newblock Differential equations driven by discontinuous paths.
\newblock Ph.D. Thesis, Imperial College of London, 1998.

\bibitem{will2}
D.~R.~E.~Williams.
\newblock Path-wise solutions of stochastic differential equations driven by 
L\'evy processes.
\newblock {\em Rev. Mat. Iberoam.} {\bf 17} (2), pp. 295-329, 2001.

\end{thebibliography}
\end{document}